\documentclass[12pt]{amsart}
\usepackage{amsfonts,amssymb,amsmath,amscd,amsthm} 
\usepackage{graphicx} 
\usepackage{stmaryrd}
\usepackage[all]{xy}

\usepackage[colorlinks,allcolors=blue]{hyperref} 

\newtheorem{thm}{Theorem}
\newtheorem{lem}[thm]{Lemma}

\newtheorem{claim}[thm]{Claim}

\theoremstyle{definition}
\newtheorem{dfn}[thm]{Definition}
\newtheorem{ex}[thm]{Example}
\newtheorem{rmk}[thm]{Remark}

\numberwithin{thm}{section}
\numberwithin{equation}{section}

\newcommand{\Proof}{\noindent {\it Proof}.\ \ }
\newcommand{\Hom}{\operatorname{Hom}}

\newcommand{\Ext}{\operatorname{Ext}}
\newcommand{\Spec}{\operatorname{Spec}}

\newcommand{\ob}{\operatorname{ob}}
\newcommand{\Hilb}{\operatorname{Hilb}}
\newcommand{\HF}{\operatorname{HF}}
\newcommand{\red}{\operatorname{red}}
\newcommand{\im}{\operatorname{im}}

\newcommand{\Bs}{\operatorname{Bs}}

\newcommand{\Pic}{\operatorname{Pic}}

\newcommand{\coker}{\operatorname{coker}}
\newcommand{\Gr}{\operatorname{Gr}}

\newcommand{\Aut}{\operatorname{Aut}}

\newcommand{\codim}[1]{\operatorname{codim}}

\renewcommand{\labelenumi}{{\rm (\arabic{enumi})}}

\newcommand{\mapright}[1]{%
  \smash{\mathop{%
      \hbox to 1cm{\rightarrowfill}}\limits^{#1}}}
\newcommand{\mapleft}[1]{%
  \smash{\mathop{%
      \hbox to 1cm{\leftarrowfill}}\limits^{#1}}}
\newcommand{\mapdown}[1]{\Big\downarrow
  \llap{$\vcenter{\hbox{$\scriptstyle#1\,$}}$ }}
\newcommand{\mapup}[1]{\Big\uparrow
  \rlap{$\vcenter{\hbox{$\scriptstyle#1$}}$ }}

\title[Obstructions to deforming space curves...]{
  Obstructions to deforming space curves
  lying on a del Pezzo surface}
\address{
  Department of Mathematical Sciences,
  Tokai University,
  4-1-1 Kitakaname, Hiratsuka, 
  Kanagawa 259-1292, JAPAN}
\author{Hirokazu Nasu}
\email{nasu@tokai.ac.jp}

\begin{document}

\begin{abstract}
We study the deformations of space curves $C \subset \mathbb P^4$,
assuming that they are contained in a smooth complete intersection
$S_{2,2} \subset \mathbb P^4$,
i.e.,~a smooth del Pezzo surface of degree $4$.
We give sufficient conditions for $C$ to be (un)obstructed
in terms of the degree $d$ and the genus $g$ of $C$.
We prove that if $d>8$, $g\ge 2d-12$, and $h^1(C,\mathcal I_C(2))=1$,
then $C$ is {\em obstructed} and {\em stably degenerate},~i.e., 
$C$ has some first order infinitesimal deformations in $\mathbb P^4$
not contained in any deformations of $S_{2,2}$ in $\mathbb P^4$, 
but they do not lift to any global deformations.
(As a result, every global deformation of $C$ in $\mathbb P^4$ 
is contained in a deformation of $S_{2,2}$ in $\mathbb P^4$.)
As an application, we construct infinitely many examples of
irreducible components of the Hilbert scheme $\Hilb^{sc} \mathbb P^4$
of smooth connected curves in $\mathbb P^4$,
along which $\Hilb^{sc} \mathbb P^4$ is generically non-reduced.
In the case $d=14$ and $g=16$,
we obtain a non-reduced component of $\Hilb^{sc} \mathbb P^4$
of dimension $55$ with $\dim T_{\Hilb^{sc} \mathbb P^4}=57$, 
analogous to Mumford's example of a non-reduced component of
$\Hilb^{sc} \mathbb P^3$, whose general member is contained in
a smooth cubic surface $S_3 \subset \mathbb P^3$.
\end{abstract}

\maketitle

\section{Introduction}
\label{sec:intro}

Let $X$ be a projective variety of dimension at least $3$,
$C$ a smooth connected curve in $X$, and
$\Hilb^{sc} X$ the Hilbert scheme of smooth connected curves in $X$.
Assuming the presence of a smooth surface $S$
such that $C \subset S\subset X$ and 
a $(-1)$-curve $E \simeq \mathbb P^1$ on $S$,
the (embedded) deformations of $C$ in a smooth $3$-fold $X$
were studied by Mukai and Nasu (\cite{Mukai-Nasu}).
The existence of $S$ and $E$ enabled them to compute the obstructions 
to deforming $C$ in $X$, and produce a number of examples of 
generically non-reduced components of $\Hilb^{sc} X$ of the $3$-folds $X$.
Later, their result has been generalized to the case where
$X$ is e.g., an (Enriques-)Fano $3$-fold,
$S$ is e.g., a del Pezzo, $K3$ or Enriques surface,
and $E$ is e.g., 
a $(-1)$-curve, a $(-2)$-curve, an elliptic curve or a half pencil
(cf.~\cite{Nasu4, Nasu5, Nasu6, Nasu7, Nasu8}).
On the other hand, in \cite{Kleppe87, Nasu1, Nasu9}, 
the deformations of space curves $C \subset \mathbb P^3$
have been studied when $C$ are contained 
in a smooth cubic surface $S_3 \subset \mathbb P^3$.
A systematic study of their deformations 
(using the coordinate of $\Pic S_3 \simeq \mathbb Z^7$)
was started by Kleppe in \cite{Kleppe87}.
By virtue of a beautiful geometry of lines on cubic surfaces,
their deformations have been well understood.
For example, 
the dimension of $\Hilb^{sc} \mathbb P^3$ at the point $[C]$
has been computed for curves $C \subset S_3 \subset \mathbb P^3$
of degree $d>9$ and genus $g \ge 3d-18$
when $C$ is $3$-normal 
(in \cite{Kleppe87}) or $2$-normal (in \cite{Nasu9}).
Here a projective variety $Y \subset \mathbb P^n$
is called {\em $m$-normal} if $H^1(\mathbb P^n,\mathcal I_Y(m))=0$,
where $\mathcal I_Y$ is the sheaf of ideals defining $Y$ in $\mathbb P^n$
and $\mathcal I_Y(m) = \mathcal I_Y \otimes_{\mathbb P^n} \mathcal O_{\mathbb P^n}(m)$.

The purpose of this paper is twofold.
The first purpose is to further develop the theories and techniques
developed in \cite{Nasu1,Mukai-Nasu,Nasu6}
and compute obstructions to deforming curves in higher 
dimensional algebraic varieties.
The second purpose is to understand the deformations of space curves
lying on a smooth del Pezzo surface of degree greater than three.
In this paper,
we mainly study the deformations of $C \subset \mathbb P^4$
contained in a smooth complete intersection $S=S_{2,2} \subset \mathbb P^4$
of two quadrics,~i.e.,~a smooth quartic del Pezzo surface in $\mathbb P^4$.
We first give a sufficient condition for such a curve to be 
{\em unobstructed} in $\mathbb P^4$, in other words,
$\Hilb^{sc} \mathbb P^4$ is nonsingular at the corresponding point.

\begin{figure}[h]
  \label{fig:d-g-region}
  \begin{center}
    \includegraphics[clip,scale=1.2]{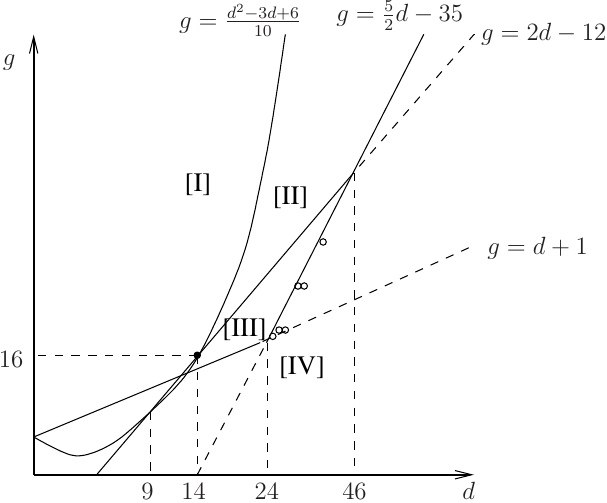}
  \end{center}
  \caption{region of $(d,g)$}
\end{figure}

\begin{thm}
  \label{thm:main1}
  Let $C \subset \mathbb P^4$ be a smooth connected curve
  of degree $d>8$ and genus $g$ contained in
  a smooth quartic del Pezzo surface $S=S_{2,2} \subset \mathbb P^4$.
  Then
  \begin{enumerate}
    \item Such curves $C$ are parametrised by
    a finite union of locally closed 
    irreducible subsets $W \subset \Hilb^{sc} \mathbb P^4$
    of the same dimension $d+g+25$.
    \item If $H^1(\mathcal I_C(2))=0$ for a member $C$ of $W$,
    then the closure $\overline W$ of $W$ in $\Hilb^{sc} \mathbb P^4$
    is an irreducible component of $\Hilb^{sc} \mathbb P^4$ and 
    $\Hilb^{sc} \mathbb P^4$ is smooth along $W$.
    \item If $H^1(\mathcal O_C(2))=0$ for a member $C$ of $W$,
    then $\Hilb^{sc} \mathbb P^4$ is smooth along $W$ and 
    $\overline W$ is of codimension $2h^1(\mathcal I_C(2))$ in $\Hilb^{sc} \mathbb P^4$.
    \item We have $H^1(\mathcal I_C(2))=0$ if $g> (d^2-3d+6)/10$ 
    (i.e.,~$(d,g)$ belongs to [I] in Figure~\ref{fig:d-g-region}),
    and $H^1(\mathcal O_C(2))=0$ if
    \begin{enumerate}
      \item $g< d+1$, or
      \item $g< 5d/2-35$ and
      \begin{equation}
	\label{eqn:exceptional pairs}
	(d,g)\not\in P:=
	\left\{(25,27),(26,28),(27,28),(29,36),(30,36),(33,45)\right\}
      \end{equation}
      (i.e.,~$(d,g)$ belongs to [IV] in Figure~\ref{fig:d-g-region}).
    \end{enumerate}
  \end{enumerate}
\end{thm}

In Theorem~\ref{thm:main1},
every member $C$ of $W$ is contained in unique $S$ by $d>8$
(cf.~Lemma~\ref{lem:uniqueness of pencil}).
The number $d+g+25$ is equal to the dimension of
the Hilbert-flag schemes $\HF \mathbb P^4$ of $\mathbb P^4$
at the nonsingular point $(C,S)$
(cf.~Lemma~\ref{lem:nonsingular of expected dimension}).
The two numbers
$h^1(\mathcal I_C(2))$ and $h^1(\mathcal O_C(2))$ are closely related to
the linear system $|C+mK_S|$ ($m=2,3$) of canonical adjunctions on $S$
by Lemma~\ref{lem:normality and speciality}.
Due to the existence of
the {\em standard coordinate} of $[C]$ in $\Pic S \simeq \mathbb Z^6$,
(cf.~Definition~\ref{dfn:standard coordinates}),
there exists a $1$-to-$1$ correspondence
between the set of all maximal irreducible families 
in $\Hilb^{sc} \mathbb P^4$
of curves of degree $d$ and genus $g$
contained in a smooth quartic del Pezzo surface and
the $6$-tuples $(a;b_1,\dots,b_5)$ of integers 
satisfying a set of numerical 
conditions \eqref{standard-prescribed-nef-big}
(cf.~Lemma~\ref{lem:maximal families and 6-tuples}).
A similar correspondence was given in \cite{Kleppe87} for
families of space curves in $\mathbb P^3$ lying on a smooth cubic surface.
Thus we can explicitly compute the numbers
$h^1(\mathcal I_C(2))$ and $h^1(\mathcal O_C(2))$
by using the coordinates $(a;b_1,\dots,b_5)\in \mathbb Z^6$ 
of $[C]$ in $\Pic S\simeq \mathbb Z^6$.
By Theorem~\ref{thm:main1}, 
if either $H^1(\mathcal I_C(2))$ or $H^1(\mathcal O_C(2))$ vanishes,
then $C$ is unobstructed and
we are able to determine the dimension of $\Hilb^{sc} \mathbb P^4$ at $[C]$.
Otherwise, $(d,g)$ belongs to the region [II] or [III] of
Figure~\ref{fig:d-g-region} except for six exceptional pairs in $P$.
In [III], by dimension reason,
there exists a global deformation of $C\subset \mathbb P^4$ 
not contained in any deformation of $S \subset \mathbb P^4$ 
(cf.~Remark~\ref{rmk:dimension reason}).
On the other hand, 
we expect that if $(d,g)$ belongs to [II] and $H^1(\mathcal I_C(2))\ne 0$,
then the maximal family $\overline W$
corresponds to a generically non-reduced component of
$\Hilb^{sc} \mathbb P^4$.
We prove that this is true if $H^1(\mathcal I_C(2))$
is of dimension one.

\begin{thm}
  \label{thm:main2}
  Let $W$ be a maximal irreducible family of smooth connected curves
  $C \subset \mathbb P^4$ of degree $d$ and genus $g$
  contained in a smooth quartic del Pezzo surface
  $S=S_{2,2} \subset \mathbb P^4$.
  If $d>8$, $g\ge 2d-12$ and $h^1(\mathcal I_C(2))=1$,
  then $\overline W$ is a component of
  $(\Hilb^{sc} \mathbb P^4)_{\red}$ and 
  $\Hilb^{sc} \mathbb P^4$ is generically non-reduced along $\overline W$.
\end{thm}

Similar results can be found in \cite{Kleppe87,Nasu1,Nasu9}
for space curves in $\mathbb P^3$ lying on smooth cubic surfaces.
The following example is an analogue of Mumford's example 
(cf.~\cite{Mumford}), which is a celebrated example 
of a non-reduced component of $\Hilb^{sc} \mathbb P^3$.

\begin{ex}
  \label{ex:main}
  Let $S:=S_{2,2} \subset \mathbb P^4$ be as above
  and $E$ a line on $S$.
  We consider a complete linear system 
  $\Lambda:=|-3K_S+2E|\simeq \mathbb P^{29}$ on $S$.
  Then $\Lambda$ contains a smooth connected member, that is 
  a curve $C \subset \mathbb P^4$ of degree $14$ and genus $16$. 
  Such curves $C$ (with variable $S$ and $E$) are parametrised by
  a locally closed irreducible subset $W$ of $\Hilb^{sc} \mathbb P^4$
  of dimension $55$ ($=14+16+25$) by Theorem~\ref{thm:main1}.
  It follows from
  the exact sequence
  $$
  0 \longrightarrow \underbrace{H^0(N_{C/S})}_{\simeq k^{29}}
  \longrightarrow H^0(N_{C/\mathbb P^4})
  \longrightarrow 
  \underbrace{H^0(N_{S/\mathbb P^4}\big{\vert}_C)}_{\simeq H^0(\mathcal O_C(2)^{\oplus 2})\simeq k^{28}}
  \longrightarrow 0
  $$
  on $C$ that $h^0(N_{C/\mathbb P^4})=57$.
  Since $h^1(\mathcal I_C(2))=1$, Theorem~\ref{thm:main2} shows that
  $\overline W$ is an irreducible component of 
  $(\Hilb^{sc} \mathbb P^4)_{\red}$ and 
  $\Hilb^{sc} \mathbb P^4$ is generically non-reduced 
  along $\overline W$.
\end{ex}
We note that this example attains the minimal degree 
and the minimal genus of curves $C \subset \mathbb P^4$ in 
the region [II].
Our proof of Theorem~\ref{thm:main2} depends on the computations of 
obstructions.
If $\tilde C \subset \mathbb P^4 \times \Spec k[t]/(t^2)$ is 
a first order infinitesimal deformation of $C \subset \mathbb P^4$, 
then there exists a global section $\alpha$ of $N_{C/\mathbb P^4}$ 
corresponding to $\tilde C$ and an element
$\ob(\alpha) \in H^1(N_{C/\mathbb P^4})$
such that $\ob(\alpha)=0$ if and only if 
$\tilde C$ is liftable to a second order deformation.
Here $\ob(\alpha)$ is called the {\em primary obstruction} of $\alpha$
and it is defined as a cup product
$$\ob:H^0(C,N_{C/\mathbb P^4}) \longrightarrow H^1(C,N_{C/\mathbb P^4}),
\quad
\alpha \longmapsto \alpha^2(=\alpha \cup \alpha).
$$
(See \S\ref{subsec:primary and exterior} for details.)
It is generally difficult to compute $\ob(\alpha)$ directly.
In this paper, we apply a theory of 
{\em exterior components, infinitesimal deformations with poles}
in \cite{Mukai-Nasu,Nasu4,Nasu5} to the computation, and
finally, reduce the computation of the primary obstruction $\ob(\alpha)$
to that of the Serre duality pairing on 
a $(-1)$-curve $E \simeq \mathbb P^1$ on $S_{2,2}$ determined by $C$.
We deduce the nonzero of $\ob(\alpha)$ from the perfectness of the pairing,
and prove $\ob(\alpha)\ne 0$
if $C$ is general in $W$ and 
$\alpha$ is not contained in the image of the tangent map
$$
p_1: T_{\HF^{sc} \mathbb P^4,(C,S)}
\longrightarrow T_{\Hilb^{sc} \mathbb P^4,[C]} =H^0(C,N_{C/\mathbb P^4})
$$
of a natural morphism
$pr_1: \HF^{sc} \mathbb P^4 \rightarrow \Hilb^{sc} \mathbb P^4$
from the Hilbert flag scheme $\HF^{sc} \mathbb P^4$.
Here $\HF \mathbb P^4$ parametrizes
all pairs $(C,S)$ of closed subschemes of $\mathbb P^4$ 
such that $C \subset S$.
This implies that $\dim \Hilb^{sc}_{[C]} \mathbb P^4
=\dim_{[C,S]} \HF^{sc} \mathbb P^4=d+g+25$ by Lemma~\ref{lem:maximality}
and $\overline W$ is maximal in $(\Hilb^{sc} \mathbb P^4)_{\red}$.
Then the generic non-reducedness of $\Hilb^{sc} \mathbb P^4$
along $\overline W$ naturally follows from the existence of
nonzero obstructions and thereby we obtain Theorem~\ref{thm:main2}.

The organization of this paper is as follows. 
In \S\ref{sec:preliminaries}, we recall some basic results 
on Hilbert-flag schemes and primary obstructions.
In \S\ref{sec:unobstructed}, we discuss the unobstructedness of
space curves lying on a smooth quartic del Pezzo surface 
and prove Theorem~\ref{thm:main1}.
In \S\ref{sec:non-reducedness} we compute obstructions to deforming
such curves and prove Theorem~\ref{thm:main2}.

\vskip 3mm

\paragraph{\bf Acknowledgment}

I would like to thank Mr.~Akira Yamada for many helpful discussions.
He carefully read a draft version of this paper and wrote a programming code
to obtain the boundaries of the regions 
[I] and [IV] in Theorem~\ref{thm:main1}, which was very helpful.
This work was supported in part by JSPS KAKENHI Grant Numbers JP24K06677.

\paragraph{\bf Notation and Conventions}

Throughout the paper, we work over an algebraically closed 
field $k$ of characteristic zero.
Given a closed subscheme $Y \subset X$,
$\mathcal I_{Y/X}$ (or $\mathcal I_Y$)
denotes the sheaf of ideals defining $Y$ in $X$,
and $N_{Y/X}$ denotes the normal sheaf of $Y$ in $X$.
Given a divisor $D$ and a sheaf $\mathcal F$ on $X$, 
$\chi(X,D)$ and $\chi(X,\mathcal F)$ represent the 
Euler characteristics of $\mathcal O_X(D)$ and $\mathcal F$, 
respectively.
A divisor $D$ on a surface $S$ is said to be nef,
if $D.C \ge 0$ for all curves $C$ on $S$, and $D$ is called big if $D^2>0$.
We write $D \ge 0$ for $D$ effective and $D>0$ for $D$ effective and nonzero.
We denote the Hilbert scheme of $X$ by $\Hilb X$ and
its open subscheme parametrizing smooth connected curves in $X$ 
by $\Hilb^{sc} X$.

\section{Preliminaries}
\label{sec:preliminaries}

\subsection{Curves on del Pezzo surfaces}
\label{subsec:delpezzo}

In this section, we recall some results on 
linear systems on del Pezzo surfaces.
We refer to e.g.~\cite{Manin,Nasu4} for proofs and
more details of results in this section.

A smooth projective surface $S$ is called {\em del Pezzo}
if its anticanonical class $-K_S$ is ample.
It is known that every del Pezzo surface is isomorphic to 
a blow-up $S_n$ of $\mathbb P^2$ at $9-n$ points ($0 < n\le 9$)
in general positions, or $\mathbb P^1\times \mathbb P^1$ 
(cf.~\cite[\S24]{Manin}).
The {\em degree} of $S$, denoted by $\deg S$, is defined to be
the self-intersection number $K_S^2$, and
$\deg S=n$ if $S\simeq S_n$ and $\deg S=8$ otherwise.
If $\deg S\ge 3$, then there exists a canonical embedding of $S$ 
into $\mathbb P^n$ ($\simeq |-K_S|^{\vee}$),
and the image $S'$ of $S$ in $\mathbb P^n$
is known to be projectively normal, i.e., $S'$ is $l$-normal 
for every integer $l \in \mathbb Z$.
A {\em line} on $S$ is a curve $\ell \simeq \mathbb P^1$ on $S$ 
with $\ell^2=-1$. 
We note that there exist no lines on $\mathbb P^2$ ($=S_9$) 
nor $\mathbb P^1\times \mathbb P^1$,
because the definition of lines implies $-K_S.\ell =1$.

We collect some basic results on divisors on del Pezzo surfaces.
Let $D$ be a divisor on a del Pezzo surface $S$ of degree $n$.
\begin{enumerate}
  \renewcommand{\labelenumi}{{\rm ({\bf P\arabic{enumi}})}}
  \item If $n \ge 8$ and $D\ge 0$ then $D$ is nef.
  If $n \le 7$ then $D$ is nef if and only if 
  $D.\ell \ge 0$ for all lines $\ell$ on $S$.
  \item If $D$ is nef then [i] $D \ge 0$, [ii] $D^2 \ge 0$,
  [iii] $|D|$ is base point free, except for the case $n=1$ and $D\sim -K_S$,
  and [iv] $D^2=0$ if and only if $D=mq$ for some $m\ge 0$ 
  and a conic $q$ on $S$.
  Here a {\em conic} means a curve $q$ on $S$ with $-K_S.q=2$ and $q^2=0$.
  \item If $D \ge 0$ then there exists a decomposition
  $$
  |D| = |D'| + F, \qquad F:=-\sum_{D.\ell<0} (D.\ell) \ell
  $$
  of the linear system $|D|$ spanned by $D$ (i.e.~a {\em Zariski decomposition}),
  where $D'$ is nef and the sum in the equation of $F$
  is taken over all (disjoint) lines $\ell$ on $S$ such that $D.\ell<0$
  (cf.~\cite[Lemma~2.2]{Nasu4}).
  We denote the fixed part $F$ of $|D|$ by $\Bs|D|$.
  \item If $D$ is nef and $\chi(S,-D)\ge 0$, then $H^1(S,-D)=0$
  (cf.~\cite[Lemma~2.1]{Nasu4}).
  If $D \ge 0$ and $H^1(S,-D)\ne 0$, then
  \begin{enumerate}
    \item $D$ is not nef, or 
    \item there exists a conic $q$ on $S$ and an integer $m\ge 2$
    such that $D \sim mq$.
  \end{enumerate}
  Then $h^1(S,-D)=m-1$ in [ii], while in [i], 
  if moreover $D^2>0$ then $h^1(S,-D)=h^0(F,\mathcal O_F)$,
  where $F$ is the fixed part of $|D|$
  (cf.~\cite[Lemma~2.4]{Nasu1}).

  \item If $D$ is nef and nonzero,
  then $|D|$ contains an irreducible curve $C$ on $S$, whose
  degree $d$ and genus $g$ satisfies
  \begin{equation}
    \label{ineq:Hodge index}
    (-K_S)^2 C^2-(-K_S.C)^2=n(d+2g-2)-d^2\le 0
  \end{equation}
  by Hodge index theorem.

  \item Suppose that $\chi(S,-D)\ge 0$ and $\deg S\ge 3$.
  If $D+K_S \ge 0$ then $D$ is big (cf.~\cite[Lemma~2.8]{Nasu9}).
  If $D\ge 0$ and $D.\ell \ge -1$ for every line $\ell$ on $S$,
  then $H^1(S,3F-D)=0$ for $F=\Bs|D|$ (cf.~\cite[Lemma~2.7]{Nasu9}).
\end{enumerate}

The normality and speciality of
space curves lying on del Pezzo surfaces
can be computed by the linear system of canonical adjunctions
on the surfaces.

\begin{lem}
  \label{lem:normality and speciality}
  Let $S \subset \mathbb P^n$ be 
  an (anti-canonically embedded) del Pezzo surface of degree $n \ge 3$,
  and $C$ a smooth connected curve on $S$.
  Put $D_m:=C+mK_S$ in $\Pic S$.
  Then there exist isomorphisms
  \begin{align}
    \label{isom:normality}
    H^1(\mathbb P^n,\mathcal I_C(m))
    & \simeq H^1(S,-D_m), \qquad \mbox{for all $m \in \mathbb Z$, and} \\
    \label{isom:speciality}
    H^1(C,\mathcal O_C(m))
    & \simeq H^2(S,-D_m)
    \simeq H^0(S,D_{m+1})^{\vee}, \qquad 
    \mbox{if $m \in \mathbb Z_{\ge 0}$}.
  \end{align}
\end{lem}
\Proof Let $m$ be any integer. There exists an exact sequence
\begin{equation}
  \label{ses:ideal sheaves}
  0 \longrightarrow \mathcal I_S(m)
  \longrightarrow \mathcal I_C(m)
  \longrightarrow \mathcal I_{C/S}(m)
  \longrightarrow 0
\end{equation}
on $S$. Since $S$ is projectively normal and 
$\mathcal I_{C/S}(m) \simeq \mathcal O_S(-D_m)$,
for the proof of \eqref{isom:normality}, it suffices to show that 
$H^2(\mathcal I_S(m))=0$ for all $m$.
Then this follows from the exact sequence
$0 \rightarrow \mathcal I_S(m) \rightarrow \mathcal O_{\mathbb P^n}(m)
\rightarrow \mathcal O_S(m) \rightarrow 0$ on 
$\mathbb P^n$ and that $H^2(\mathcal O_{\mathbb P^n}(m))=H^1(\mathcal O_S(m))=0$.
Since $H^2(\mathcal O_S(m))=0$ for $m \ge 0$,
the first isomorphism in \eqref{isom:speciality} is deduced from the exact sequence
$$
0 \longrightarrow \mathcal I_{C/S}(m)
\longrightarrow \mathcal O_S(m)
\longrightarrow \mathcal O_C(m)
\longrightarrow 0
$$ on $S$, while the last isomorphism in \eqref{isom:speciality} follows from
the Serre duality.
\qed

\medskip

Let $D_m$ denote the class of $C+mK_S$ in $\Pic S$.
Then provided that $D_m\ge 0$ and $D_m^2>0$,
$C$ is $m$-normal if and only if $D_m$ is nef, otherwise
$$
h^1(\mathbb P^n,\mathcal I_C(m))=h^0(F,\mathcal O_F),
$$
where $F$ is the fixed part of $|D_m|$.

We next setup some notations concerning the
coordinates of divisors on del Pezzo surfaces.
For the sake of simplicity, in what follows we assume that $\deg S=4$, 
i.e., $S$ is isomorphic to
a blow up of $\mathbb P^2$ at $5$ points in a general position.
Then the Picard group $\Pic S$ of $S$ is a free abelian group
generated by the class $\mathbf l$ of the pullback of lines in $\mathbb P^2$
and $5$ exceptional curves $\mathbf e_1, \dots, \mathbf e_5$ on $S$.
Thus every divisor $D \sim a \mathbf l - \sum_{i=1}^5 b_i \mathbf e_i$ on $S$
corresponds to a $6$-tuple $(a;b_1,\ldots,b_5) \in \mathbb Z^6$ of integers
by coordinates.
The anticanonical class $-K_S$ corresponds to $(3;1,\dots,1)$.
Then we take actions of the Weyl group into account.
There exists a Weyl group $W \subset \Aut(\Pic S)$
corresponding to the root system of Dynkin type $D_4$,
which is generated by the permutations of $\mathbf e_i$ ($1 \le i \le 5$)
and the Cremona transformation $\sigma$ on $\mathbb P^2$.
Then by the action of the Weyl group,
there exists a suitable blow-up $S \rightarrow \mathbb P^2$ such that
\begin{equation}\label{eqn:standard}
  b_1 \ge \dots \ge b_5 \quad \mbox{and} \quad a \ge b_1 + b_2 + b_3
\end{equation}
(see e.g.\cite{Manin},\cite[\S5.3]{Nasu4},\cite{Nasu9}).
\begin{dfn}
  \label{dfn:standard coordinates}
  We say the basis
  $\left\{\mathbf l,\mathbf e_1,\ldots,\mathbf e_5\right\}$ of $\Pic S$
  (or the coordinate $(a;b_1,\dots,b_5)$ of $D$ in $\Pic S$) is
  {\em standard} (for $D$) if \eqref{eqn:standard} is satisfied.
\end{dfn}

Under the standard coordinate, $D$ is nef if and only if $b_5\ge 0$.
If $D$ is nonzero and nef, then $|D|$ contains an irreducible curve as a member,
whose degree $d$ and (arithmetic) genus $g$ are respectively obtained 
by the formulas
\begin{equation}\label{eqn:degree and genus}
d = 3a- \sum_{i=1}^5 b_i \quad \mbox{and} \quad 
g = \dfrac{(a-1)(a-2)}2 - \sum_{i=1}^5 \dfrac{b_i(b_i-1)}2.
\end{equation}
We have by Hodge index theorem that $0 \le g \le 1+(d-4)d/8$.

\subsection{Hilbert-flag schemes and maximal families}
\label{subsec:flag and max_fam}

In this section, we recall basic results on Hilbert-flag schemes and 
(primary) obstructions associated with first order infinitesimal 
deformations of closed subschemes. Given a projective scheme $X$,
the {\em Hilbert-flag scheme} $\HF X$ of $X$ parametrizes all pairs $(C,S)$ 
of closed subschemes of $X$ such that $C \subset S \subset X$
(cf.~\cite{Kleppe87,Sernesi}).
$\HF X$ represents the Hilbert-flag functor, which
assigns to each base scheme $B$
a flat family $\mathcal C \subset \mathcal S \subset X\times_k B$
over $B$, whose fibers are pairs $(C,S)$ of closed subschemes such that
$C \subset S \subset X$ with fixed Hilbert polynomials.
There exist two natural morphisms $pr_i: \HF X \rightarrow \Hilb X$ ($i=1,2$)
corresponding to the two projections 
$(C,S) \mapsto [C]$ and $(C,S) \mapsto [S]$, respectively.

The normal sheaf $N_{(C,S)/X}$ of $(C,S)$ in $X$ is defined by
the fiber product of $N_{C/X}$ and $N_{S/X}$ over $N_{S/X}\big{\vert}_C$.
In what follows, for simplicity, we assume that the closed embeddings
$C \hookrightarrow S$ and $S\hookrightarrow X$ are both regular.
Then:
\begin{enumerate}
  \renewcommand{\labelenumi}{{\rm ({\bf P\arabic{enumi}})}}
  \setcounter{enumi}{6}
  \item The cohomology groups $H^i(N_{(C,S)/X})$ ($i=0,1$) represent
  the tangent space and the obstruction space of $\HF X$ at $(C,S)$,
  respectively (cf.~\cite[Proposition~4.5.3]{Sernesi}).
  \item If $H^1(N_{(C,S)/X})=0$ then 
  $\HF X$ is nonsingular at $(C,S)$ of expected dimension,
  which coincides with $\chi(N_{(C,S)/X})$
  provided that $H^i(N_{(C,S)/X})=0$ for all $i \ge 2$.
  \item There exist two natural exact sequences
  \begin{align}
    \label{ses:normal1}
    0 \longrightarrow \mathcal I_{C/S}\otimes_S N_{S/X} 
    \longrightarrow & N_{(C,S)/X} \overset{\pi_1}{\longrightarrow}
    N_{C/X} \longrightarrow 0\\
    \label{ses:normal2}
    0 \longrightarrow N_{C/S} \longrightarrow & N_{(C,S)/X} 
    \overset{\pi_2}{\longrightarrow} N_{S/X} \longrightarrow 0
  \end{align}
  on $X$, where $\pi_i$ ($i=1,2$)
  induce the tangent maps (and the maps on obstruction spaces) of $pr_i$,
  respectively (cf.~\cite[\S2.2]{Nasu6}).
  \item If $H^1(\mathcal I_{C/S}\otimes_S N_{S/X})=0$,
  then $pr_1$ is smooth at $(C,S)$, while
  so is $pr_2$ if $H^1(N_{C/S})=0$. 
  (cf.~\cite[Lemma~A10]{Kleppe87}, \cite[Theorem~1.3.4]{Kleppe81}).
  \item Let $\mathcal I_{C/S}N_S$ denote the sheaf $\mathcal I_{C/S}\otimes_S N_{S/X}$ on $S$.
  If $H^0(\mathcal I_{C/S}N_S)=0$
  and $H^i(N_{(C,S)/X})=0$ for $i=1,2$,
  then 
  \begin{equation}
    \label{ineq:codimension of HF in Hilb}
    h^1(\mathcal I_{C/S}N_S)-h^2(\mathcal I_{C/S}N_S)
    \le 
    \dim_{[C]} \Hilb X - \dim_{(C,S)} \HF X
    \le
    h^1(\mathcal I_{C/S}N_S),
  \end{equation}
  where the inequality to the right is strict if and only if
  $\Hilb X$ is singular at $[C]$ (cf.~\cite[Theorem~2.4]{Nasu6}).
  
\end{enumerate}

The following lemma is a generalization of \cite[Lemma~2.15]{Nasu9}.

\begin{lem}[{cf.~\cite{Kleppe87},\cite{Nasu9}}]
  \label{lem:nonsingular of expected dimension}
  Let $S \subset \mathbb P^n$ be 
  an (anti-canonically embedded) del Pezzo surface of degree $n \ge 3$,
  and $C$ a smooth connected curve on $S$ of degree $d$ and genus $g$.
  Then $\HF \mathbb P^n$ is nonsingular at $(C,S)$ of expected dimension
  $$\chi(\mathbb P^n,N_{(C,S)/\mathbb P^n})=d+g+n^2+9.$$
  In fact, all the higher cohomology groups of 
  $N_{C/S},N_{S/\mathbb P^n}$ (and hence $N_{(C,S)/\mathbb P^n}$) vanish.
\end{lem}
\Proof
By adjunction, we have $N_{C/S} \simeq -K_S\big\vert_C+K_C$,
where $-K_S\simeq \mathcal O_S(1)$ is ample.
Therefore $H^i(N_{C/S})=0$ for all $i>0$ and
$\chi(N_{C/S})=d+g-1$ by Riemann-Roch theorem on $C$.
There exists an exact sequence
$$
0 \longrightarrow \mathcal O_S
\longrightarrow \mathcal O_S(1)^{\oplus n+1}
\longrightarrow T_{\mathbb P^n}\big{\vert}_S
\longrightarrow 0
$$
on $S$, which is the restriction to $S$ of the Euler sequence on $\mathbb P^n$.
Since all the higher cohomology groups of $\mathcal O_S$ and $\mathcal O_S(1)$ vanish,
so do those of $T_{\mathbb P^n}\big{\vert}_S$.
Then by additivity, $\chi(T_{\mathbb P^n}\big{\vert}_S)=
(n+1)\chi(\mathcal O_S(1))-\chi(\mathcal O_S)
=(n+1)^2-1=n^2+2n$. Since $S$ is smooth and Fano, we have
$H^i(T_S)=0$ for $i\ge 2$.
Then it follows from the exact sequence
$$
0 \longrightarrow T_S \longrightarrow T_{\mathbb P^n}\big{\vert}_S 
\longrightarrow N_{S/\mathbb P^n}
\longrightarrow 0
$$
that $H^i(N_{S/\mathbb P^n})=0$ for all $i>0$ and
$\chi(N_{S/\mathbb P^n})
=\chi(T_{\mathbb P^n}\big{\vert}_S)-\chi(T_S)
=n^2+2n-(2n-10)=n^2+10$.
Therefore by \eqref{ses:normal2},
all the higher cohomology groups of $N_{(C,S)/\mathbb P^n}$ vanish,
and we compute by additivity that
$\chi(N_{(C,S)/\mathbb P^n})=(d+g-1)+n^2+10=d+g+n^2+9$.
\qed

\medskip

We next recall the definition of $S$-maximal families of curves
(cf.~\cite{Mukai-Nasu,Nasu5,Nasu6}).
In what follows, we assume that $C$ is a smooth connected curve
and $\HF X$ is nonsingular at $(C,S)$.
Let $\Hilb^{sc} X$ denote the Hilbert scheme of smooth connected
curves in $X$ and let $pr_1': \HF^{sc} X \rightarrow \Hilb^{sc} X$ be
the restriction of $pr_1$ to $\HF^{sc} X:=pr_1^{-1} (\Hilb^{sc} X)$.
Then by assumption, there exists a unique irreducible component
$\mathcal W_{C,S}$ of $\HF^{sc} X$ passing through $(C,S)$.
We consider the image $W_{C,S}$ of $\mathcal W_{C,S}$ by $pr_1'$.

\begin{dfn}
  \label{dfn:maximal family}
  An irreducible closed subset $W_{C,S}:=pr_1'(\mathcal W_{C.S})$
  of $\Hilb^{sc} X$ is called the {\bf $S$-maximal family} 
  of curves containing $C$.
\end{dfn}

Since the smoothness is an open property,
if $pr_1$ is smooth at $(C,S)$ then so is $pr_1'$.
Then since smooth morphisms map every generic point to a generic point, 
we see that
if $H^1(X,N_{(C,S)/X})=0$ and $H^1(S,\mathcal I_{C/S}\otimes_S N_{S/X})=0$,
then the $S$-maximal family $W_{C,S}$ is a generically smooth irreducible component
of $\Hilb^{sc} X$ (cf.~\cite[Theorem~2.4]{Nasu6}).
Otherwise, one can deduce an exact sequence 
\begin{equation}
  \label{ses:first tangent map}
  \begin{CD}
    H^0(X,N_{(C,S)/X})
    \overset{p_1}{\longrightarrow} H^0(C,N_{C/X})
    \longrightarrow H^1(S,\mathcal I_{C/S}\otimes_S N_{S/X})
    \longrightarrow H^1(X,N_{(C,S)/X})
  \end{CD}
\end{equation}
from \eqref{ses:normal1}.
The following lemma is necessary for our proof of Theorem~\ref{thm:main2}.
\begin{lem}[{\cite[Lemma~2.17]{Nasu9}}]
  \label{lem:maximality}
  Suppose that $H^1(N_{(C,S)/X})=0$.
  If the primary obstruction $\ob(\alpha)$ is nonzero in $H^1(N_{C/X})$
  for every global section $\alpha \in H^0(N_{C/X})\setminus \im p_1$,
  then 
  \begin{enumerate}
    \item the $S$-maximal family $W_{C,S}$ is an irreducible 
    component of $(\Hilb^{sc} X)_{\red}$, and
    \item if $H^0(\mathcal I_{C/S}\otimes_S N_{S/X})=0$,
    then
    $$
    \dim_{[C]} \Hilb^{sc} X = \dim_{(C,S)} \HF^{sc} X.
    $$
  \end{enumerate}
\end{lem}
Here $\ob(\alpha)$ represents the primary obstruction to extend
$\alpha$ (or the first order infinitesimal deformation $\tilde C$ 
of $C$ in $X$ corresponding to $\alpha$), 
which we will recall in the next section.

\subsection{Primary obstructions and exterior components}
\label{subsec:primary and exterior}

Given a closed subscheme $C \subset X$, there exists
a natural $1$-to-$1$ correspondence 
between the set of first order infinitesimal deformations
$\tilde C \subset X \times_k k[t]/(t^2)$ of $C$ in $X$ 
and the set of global sections $\alpha$ of the normal sheaf $N_{C/X}$ of $C$.
There exists a unique element $\ob(\alpha)$ of 
$\Ext^1(\mathcal I_C,\mathcal O_C)$
for each $\alpha$, i.e.,~the {\em primary obstruction} of $\alpha$,
such that $\ob(\alpha)$ is zero if and only if
$\tilde C$ extends to a deformation 
${\tilde {\tilde C}}$ of $C$ over $k[t]/(t^3)$. 
For simplicity, in what follows, we assume $C \hookrightarrow X$ is regular.
Then there exist no local obstructions to deforming $C$ in $X$
and $\ob(\alpha)$ is contained in 
$H^1(N_{C/X}) \subset \Ext^1(\mathcal I_C,\mathcal O_C)$
(cf.~\cite[Theorem~4.3.5]{Sernesi}).
For the last inclusion, we recall that a natural surjection
$\mathcal I_C \twoheadrightarrow N_{C/X}^{\vee}$ of sheaves on $X$
induces an injective map
$$
  H^1(N_{C/X})\simeq \Ext^1_C(N_{C/X}^{\vee},\mathcal O_C)  \hookrightarrow
  \Ext^1_X(\mathcal I_C,\mathcal O_C).
$$
A similar map induces isomorphisms
$H^0(N_{C/X})\simeq \Hom_C(N_{C/X}^{\vee},\mathcal O_C) \simeq
\Hom_X(\mathcal I_C,\mathcal O_C)$
of cohomology groups, and we can identify a global section of $N_{C/X}$
with a global $\mathcal O_X$-linear homomorphism
from $\mathcal I_C$ to $\mathcal O_C$.
Similarly, we identify $\ob(\alpha)\in H^1(N_{C/X})$ with
an extension class of $\mathcal I_C$ by $\mathcal O_C$.
Then under the identifications, 
we have $\ob(\alpha)=\alpha \cup \mathbf e \cup \alpha$,
where $\mathbf e \in \Ext^1(\mathcal O_C,\mathcal I_C)$ is 
the extension class of the exact sequence
\begin{equation}
  \label{ses:standard}
  0 \rightarrow \mathcal I_C \rightarrow \mathcal O_X
  \rightarrow \mathcal O_C \rightarrow 0
\end{equation}
on $X$ (cf~\cite[Theorem~2.1]{Nasu5}). We give alternate expression 
of $\ob(\alpha)$ for later computations.
Let $\alpha': N_{C/X}^{\vee} \rightarrow \mathcal O_C$ be 
the $\mathcal O_C$-linear homomorphism corresponding to $\alpha$,
and $\mathbf e' \in \Ext^1_X(\mathcal O_C,N_{C/X}^{\vee})$
the extension class of the exact sequence
\begin{equation}
  \label{ses:conormal}
  0 \rightarrow N_{C/X}^{\vee} \rightarrow \mathcal O_X/\mathcal I_C^2
  \rightarrow \mathcal O_C \rightarrow 0
\end{equation}
on $X$. Then 
$\ob(\alpha)=\alpha' \cup \mathbf e' \cup \alpha'$
in $H^1(N_{C/X})\simeq \Ext^1_C(N_{C/X}^{\vee},\mathcal O_C)$.
We remark that $\ob(\alpha)\ne 0$ implies
that $\Hilb^{sc} X$ is singular at $[C]$
by infinitesimal lifting property of smoothness.

We next consider deformations of closed subschemes
with respect to a given fixed intermediate closed subscheme.
Let $S \subset X$ be a closed subscheme containing $C$
such that the two embeddings $C \hookrightarrow S \hookrightarrow X$ are both regular.
Then a natural projection $\pi_{C/S}: N_{C/X} \rightarrow N_{S/X}\big{\vert}_C$
induces two maps
$H^i(\pi_{C/S}): H^i(N_{C/X}) \rightarrow H^i(N_{S/X}\big{\vert}_C)$ ($i=0,1$)
on cohomology groups. Then given a global section $\alpha$ of $N_{C/X}$, 
the images of $\alpha$ and $\ob(\alpha)$ by the maps
are called their {\em exterior components} 
and denoted by $\pi_{C/S}(\alpha)$ and $\ob_S(\alpha)$, respectively.
For the proof of Theorem~\ref{thm:main2},
we deal with global sections of $N_{S/X}\big{\vert}_C$
which do not lift to those of $N_{S/X}$ but those of $N_{S/X}(E)$
for some effective Cartier divisor $E\ge 0$ on $S$.
The following definition is due to \cite{Mukai-Nasu}.
We refer to \cite{Nasu5,Nasu8} for more results
on infinitesimal deformations with poles.
\begin{dfn}
\label{dfn:infinitesimal deformation with pole}
Let $E$ be a nonzero effective divisor on $S$,~i.e., $E> 0$.
Then a rational section 
$\beta \in H^0(S,N_{S/X}(E))\setminus H^0(S,N_{S/X})$
is called an {\em infinitesimal deformation with pole}.
\end{dfn}
Let $\mathbf k_C \in \Ext^1(\mathcal O_C,\mathcal O_S(-C))$ denote
the extension class of the exact sequence
\begin{equation}
  \label{ses:CE}
  0 \longrightarrow \mathcal O_S(-C)
  \longrightarrow \mathcal O_S
  \longrightarrow \mathcal O_C
  \longrightarrow 0,
\end{equation}
and the class $\mathbf k_E$ similarly.
In what follows, given a coherent sheaf $\mathcal F$ on $S$,
an effective divisor $D\ge 0$ on $S$, an integer $i \ge 0$
and a cohomology class $*$ in $H^i(S,\mathcal F)$,
we denote by $r(*,D)$ the image of $*$ 
by the natural map  $H^i(S,\mathcal F) \rightarrow H^i(S,\mathcal F(D))$.
We use similar notation for $C$ also.
By Lemma~\ref{lem:restriction to C and E} below,
given a global section $\gamma$ of $N_{S/X}\big{\vert}_C$, 
$r(\gamma,E)$ lifts to a global section $\beta$ of $N_{S/X}(E)$
(an infinitesimal deformation with pole)
if and only if $r(\gamma,E) \cup \mathbf k_C=0$.

\begin{lem}
  \label{lem:restriction to C and E}
  Let $X$ be a projective scheme, $S \subset X$ a surface,
  and $C$ an irreducible curve such that
  $C \hookrightarrow S\hookrightarrow X$ are regular embeddings.
  Let $E>0$ be a divisor on $S$ whose support does not contain $C$,
  $\mathcal F$ a locally free sheaf on $S$ and 
  $\gamma$ a global section of $\mathcal F\big{\vert}_C$. Then
  \begin{enumerate}
    \item $r(\gamma,E)$ lifts to a section 
    $\beta \in H^0(S,\mathcal F(E))$ on $S$
    if and only if $r(\gamma,E) \cup \mathbf k_C=0$ 
    in $H^1(S,\mathcal F(E-C))$.
    \item Suppose that $r(\gamma,E)$ lifts to a global section $\beta$ of $\mathcal F(E)$,~i.e.,
    $$
    r(\gamma,E) = \beta\big{\vert}_C, \quad \mbox{in} \quad H^0(C,\mathcal F(E)\big{\vert}_C),
    $$
    and let $\beta\big{\vert}_E \in H^0(E,\mathcal F(E)\big{\vert}_E)$ be
    the principal part of $\beta$ along $E$. Then
    \begin{enumerate}
      \item  there exists a global section $\beta'_E$ of $\mathcal F(E-C)\big\vert_E$
      such that $r(\beta'_E,C)=\beta\big{\vert}_E$. In particular, 
      $\beta$ is contained in $H^0(S,\mathcal I_{C \cap E/S}\otimes \mathcal F(E))$.
      \item $\beta\big{\vert}_E$ is nonzero 
      if $\gamma \cup \mathbf k_C \ne 0$ in $H^1(S,\mathcal F(-C))$.
      \item $\beta'_E \cup \mathbf k_E=\gamma \cup \mathbf k_C$ in $H^1(S,\mathcal F(-C))$.
    \end{enumerate}
  \end{enumerate}
\end{lem}
\Proof The proof is very similar to that of \cite[Lemma~3.1]{Nasu5}.
(1) follows from the short exact sequence $\eqref{ses:CE}\otimes \mathcal F(E)$,
whose coboundary map coincides with the cup product map
$\cup \mathbf k_C: H^0(C,\mathcal F(E)\big{\vert}_C) \rightarrow H^1(S,\mathcal F(E-C))$
with $\mathbf k_C$.
For (2), [i] and [ii] follow from a diagram chase on the commutative diagram
\begin{equation}
  \label{diag:comm1}
  \begin{array}{cccccc}
    & 0 & & 0 & & 0\\
    & \downarrow & & \downarrow & & \downarrow \\
    0 \longrightarrow & H^0(S,\mathcal F(-C)) & \longrightarrow{} & H^0(S,\mathcal F(E-C)) & \mapright{|_E} & H^0(E,\mathcal F(E-C)\big{\vert}_E) \\
    & \mapdown{} & & \mapdown{} & & \mapdown{} \\
    0 \longrightarrow & H^0(S,\mathcal F)& \longrightarrow & H^0(S,\mathcal F(E)) & \mapright{|_E} & H^0(E,\mathcal F(E)\big{\vert}_E)\\
    & \mapdown{|_C} & & \mapdown{|_C} & & \mapdown{|_C} \\
    0 \longrightarrow & H^0(C,\mathcal F\big{\vert}_C) & \longrightarrow & H^0(C,\mathcal F(E)\big{\vert}_C) & \mapright{|_E} & \mathcal F \otimes_S k(C \cap E)\\
    & \mapdown{\cup \mathbf k_C} & & \mapdown{\cup \mathbf k_C} & &  \\
    & H^1(S,\mathcal F(-C))& \longrightarrow & H^1(S,\mathcal F(E-C)) &  &  \\
  \end{array} 
\end{equation}
of cohomology groups, which is exact both vertically and horizontally.
Finally [iii] follows from \cite[Lemma~2.8]{Mukai-Nasu}.
\qed

\section{Unobstructed curves and smooth components of $\Hilb^{sc} \mathbb P^4$}
\label{sec:unobstructed}

We first prepare two Lemmas~\ref{lem:2-normal} and \ref{lem:2-nonspecial},
which are necessary for our proof of Theorem~\ref{thm:main1}.
Let $C \subset \mathbb P^4$ be a smooth connected curve of degree $d$
and genus $g$. We suppose that $C$ is contained in
a smooth complete intersection $S=S_{2,2}$ of two quadrics in $\mathbb P^4$,
i.e.~a smooth del Pezzo surface $S \simeq S_4 \subset \mathbb P^4$
of degree $4$ (cf.~\cite[Prop.~IV.16]{Beauville}).
Let $(a;b_1,\dots,b_5)$ be the standard coordinate of $C$ 
in $\Pic S \simeq \mathbb Z^6$ 
(see Definition~\ref{dfn:standard coordinates})
and put $D_n:=C+nK_S$ for $n \in \mathbb Z$.
Since every element of the Weyl group preserves the class 
(and hence the coordinate) of $K_S$,
provided that $D_n \ge 0$ and $D_n^2>0$, 
$C$ is $n$-normal if and only if $b_5\ge n$
by Lemma~\ref{lem:normality and speciality}
\footnote{We note here 
that if $b_5 \ge n$, then $D_n$ is nef and hence $D_n\ge 0$.}.
In particular, $C$ is $2$-normal if $b_5 \ge 2$ and $D_2^2>0$.
We give a sufficient condition for $C$ to be $2$-normal
in terms of $d$ and $g$.
\begin{lem}
  \label{lem:2-normal}
  If $g > (d^2-3d+6)/10$ then $C$ is $2$-normal, 
  i.e., $H^1(\mathbb P^4,\mathcal I_C(2))=0$.
\end{lem}

\Proof
We may suppose that $d\ge 2$, since every line on $\mathbb P^4$ is 
projectively normal. We prove that $10g-d^2+3d-6\le 0$
if $C$ is not $2$-normal. Then by Lemma~\ref{lem:normality and speciality},
we can instead assume that $H^1(S,-D_2)\ne 0$, where $D_2:=C+2K_S$ in $\Pic S$.

We first consider the case where $D_2$ is nef. Then $D_2\ge 0$ and
there exists a conic $q$ on $S$ and an integer $m \ge 2$
such that $D_2 \sim mq$ by \S\ref{subsec:delpezzo}({\bf P4}).
Hence $C$ has the standard coordinate $(m+6;m+2,2,2,2,2)$ 
in $\Pic S\simeq \mathbb Z^6$.
We compute that $d=2m+8$ and $g=3m+5$ by \eqref{eqn:degree and genus}.
Thus we see that $10g-d^2+3d-6=4+4m-4m^2\le 0$ by $m \ge 2$.

We next consider the case where $D_2$ is not nef. 
Then there exists a line $\ell$
on $S$ such that $C.\ell=(D_2-2K_S).\ell < 2$, which implies 
that the last coordinate $b_5$ of $C$ in $\Pic S$ is either $0$ or $1$
by the standard basis of $\Pic S$.
Let $\varepsilon: S \rightarrow S'$ be the blow-down of 
the line $\mathbf e_5$, i.e., a $(-1)$-curve on $S$.
Then $S'$ is a smooth del Pezzo surface of degree $5$.
Moreover, $C$ is mapped by $\varepsilon$ 
isomorphically onto a curve $C'$ on $S'$
of degree either $d'=d+1$ or $d'=d$, as $b_5=1$ or $b_5=0$.
If $b_5=1$, then by Hodge index theorem, we obtain
$$
(-K_{S'})^2C'^2-(-K_{S'}.C')=5(d'+2g-2)-d'^2
=10g-d^2+3d-6\le 0.
$$
Similarly, if $b_5=0$ (i.e.,~$d'=d$)
then we have $10g-d^2+5d-10\le 0$ and hence
$10g-d^2+3d-6\le 4-2d \le 0$ by $d \ge 2$.
Thus we have completed the proof.
\qed

\medskip

We also give a sufficient condition for 
the invertible sheaf $\mathcal O_C(2)$ on $C$ to be non-special
in terms of $d$ and $g$ (with $a,b_1,\dots,b_5$).

\begin{lem}
  \label{lem:2-nonspecial}
  If $g<\frac{5}{2}d-35$ and $C$ does not belong to
  the $6$ exceptional classes 
  \begin{equation}
    \label{ineq:exceptional classes}
    (\lambda;\mu,0,0,0,0), 
    \qquad (9 \le \lambda < 12, 0 \le \mu < 12-\lambda),
  \end{equation}
  in $\Pic S \simeq \mathbb Z^6$ under a standard basis for $C$,
  then $\mathcal O_C(2)$ is non-special
  i.e., $H^1(C,\mathcal O_C(2))=0$.
  Moreover, if $C$ belongs to the exceptional classes 
  then $(d,g)$ belongs to the set $P$ in \eqref{eqn:exceptional pairs}.
\end{lem}

\Proof
We suppose that $\mathcal O_C(2)$ is special and prove that
the number $N:=2g-5d+70$ is greater than or equal to $0$
unless $C$ belongs to the $6$ exceptional classes
in \eqref{ineq:exceptional classes}.
By Lemma~\ref{lem:normality and speciality}, we see that
the divisor $D_3:=C+3K_S$ on $S$ is effective.
By Zariski decomposition, we have
$|D_3|=|D_3'|+F$, where $D_3'$ is nef and $F$ is the fixed part of $|D_3|$.
Then 
$$
F=\sum_{i=1}^e m_i \ell_i, \quad \mbox{and} \quad F^2=-\sum_{i=1}^e m_i^2,
$$
where the sum is taken over
all lines $\ell_i$ on $S$ such that $(0 \le) C.\ell_i< 3$
with coefficients $m_i:=-(D_3.\ell_i)$.
Since $\ell_i$'s are mutually disjoint and $S$ is a $5$-points blow-up 
of $\mathbb P^2$, the number $e$ of such lines is at most $5$. 
We compute that
$$
D_3^2=(C+3K_S)^2=C^2+6C.K_S+9K_S^2=2g-2-5d+36=2g-5d+34.
$$
It follows from $(D_3').F=0$ that
\begin{equation}
  N=D_3^2+36=(D_3')^2+F^2+36.
\end{equation}
Since $0< m_i\le 3$ for every $i$, 
if $e\le 4$ then $F^2\ge 4(-9)=-36$ and hence $N\ge 0$ by $(D_3')^2\ge 0$.

Suppose that $e=5$. Then there exists a standard basis
$\left\{\mathbf l,\mathbf e_1,\dots,\mathbf e_5\right\}$ of $\Pic S$
such that $\left\{\mathbf e_1,\dots,\mathbf e_5\right\}
=\left\{\ell_1,\dots,\ell_5\right\}$.
Then $D_3'=(a-9) \mathbf l$ with $a \ge 9$ because $D_3'\ge 0$.
Moreover we have $N=(a-9)^2+F^2+36$.
This implies that if $a \ge 12$ then $N\ge 0$, since $F^2\ge 5(-9)=-45$.
Suppose that $a<12$.
Then we still have $N\ge 0$
if there exists $2\le i\le 5$ such that $b_i\ne 0$.
Thus we may suppose that $b_i=0$ for all $2\le i\le 5$.
Then $N=(a-9)^2-(3-b_1)^2$ and $N<0$ if and only if $a-9 < 3-b_i$.
Thus $(a,b_1)=(\lambda,\mu)$ with
$9 \le \lambda < 12$ and $0 \le \mu < 12-\lambda$.
Thus the proof has been completed.
\qed

\medskip

We give an explicit description of maximal families of space curves
lying on a smooth quartic del Pezzo surface.

\begin{lem}
  \label{lem:uniqueness of pencil}
  If $d > 8$, then $C \subset \mathbb P^4$ 
  is contained in a unique smooth complete intersection 
  $S=S_{2,2} \subset \mathbb P^4$ and $H^0(S,N_{S/\mathbb P^4}(-C))=0$.
\end{lem}
\Proof
Since $(-2K_S-C).(-K_S)=8-d<0$, we see that
$-2K_S-C$ is not effective, and hence
$H^0(\mathcal I_{C/S}\otimes_S N_{S/\mathbb P^4})=0$ by
\begin{equation}
  \label{isom:normal bundle of S}
  \mathcal I_{C/S}\otimes_S N_{S/\mathbb P^4}
  \simeq N_{S/\mathbb P^4}(-C)
  \simeq \mathcal O_S(-2K_S-C)^{\oplus 2}.
\end{equation}
Then the exact sequence \eqref{ses:ideal sheaves} induces
an isomorphism between the linear system
$|\mathcal I_C(2)|$ of quadrics in $\mathbb P^4$ containing $C$ 
and the pencil $|\mathcal I_S(2)|$ 
of quadrics defining $S$ in $\mathbb P^4$,
and thus $S$ is uniquely determined by $C$.
\qed

\begin{lem}
  \label{lem:uniqueness of 6-tuple}
  Let $W \subset \Hilb^{sc} \mathbb P^4$ be a maximal irreducible family
  of smooth connected curves $C \subset \mathbb P^4$ of degree $d>8$
  contained in a (unique) smooth complete intersection 
  $S=S_{2,2} \subset \mathbb P^4$.
  Then the standard coordinates of $[C] \in \Pic S\simeq \mathbb Z^6$ 
  are same for all $C \in W$ (and $S \supset C$).
\end{lem}
\Proof
Let $C \in W$ be a member and $S \subset \mathbb P^4$ 
the smooth complete intersection $S_{2,2}$ containing $C$.
We consider the standard coordinate $(a;b_1,\dots,b_5)\in \mathbb Z^6$ of 
$[C]$ with respect to a standard basis
$\left\{\mathbf l,\mathbf e_1,\dots,\mathbf e_5\right\}$ of $\Pic S$.
Let $C'$ be any member of $W$ and $S' \subset \mathbb P^4$
the smooth complete intersection containing $C'$.
Then since $H^1(\mathcal O_S)=0$, the Picard group of $S$ does not change
under smooth deformations of $S$ and hence we have $\Pic S\simeq \Pic S'$.
Since $H^1(\mathcal O_S(\ell))=0$ for every line $\ell$ on $S$,
the classes of 5 lines $\mathbf e_1,\dots,\mathbf e_5$ on $S$
deform to those $\mathbf e_1',\dots,\mathbf e_5'$ on $S'$
such that $\mathbf e_i'.\mathbf e_j'=0$ for $i\ne j$,
and similarly for $\mathbf l=[\mathcal O_{\mathbb P^2}(1)] \in \Pic S$.
Since the intersection numbers of $C$ with the basis of $\Pic S$
are preserved under the deformations of $S$,
$C$ and $C'$ have the same (standard) coordinate in $\mathbb Z^6$.
\qed

\begin{lem}[{cf.~\cite[Remark 2]{Kleppe87}}]
  \label{lem:maximal families and 6-tuples}
  Let $d>8$ and $g$ be two integers satisfying \eqref{ineq:Hodge index} with $n=4$.
  Then the maximal irreducible families
  of smooth connected curves $C \subset \mathbb P^4$
  of degree $d$ and genus $g$ contained in 
  a smooth complete intersection $S_{2,2} \subset \mathbb P^4$
  are in $1$-to-$1$ correspondence with
  $6$-tuples $(a;b_1,\dots,b_5)$ of integers satisfying a set of conditions
  \begin{equation}
    \label{standard-prescribed-nef-big}
    \mbox{\eqref{eqn:standard}, \eqref{eqn:degree and genus},
      $b_5 \ge 0$ and $a > b_1$}.
  \end{equation}
  Moreover, every such families have the same dimension $d+g+25$.
\end{lem}
\Proof
Let $W \subset \Hilb^{sc} \mathbb P^4$ be a maximal irreducible 
family of such curves. 
Then by Lemma~\ref{lem:uniqueness of 6-tuple}, 
every member $C$ of $W$ is contained in a unique
$S=S_{2,2} \subset \mathbb P^4$ with the same
standard coordinate $(a;b_1,\dots,b_5)$ in $\Pic S$.
Since $C$ is neither a line nor a conic, we see that $C$ is nef and big.
This implies $b_5 \ge 0$ and $a>b_1$.

Conversely, let such a $6$-tuple $(a;b_1,\dots,b_5)$ of integers
be given.
Then on a smooth complete intersection $S=S_{2,2} \subset \mathbb P^4$,
the complete linear system $|D|$ spanned by the corresponding divisor
$D=a \mathbf l - \sum_{i=1}^5 b_i \mathbf e_i$ is free by $b_5 \ge 0$
(under a standard basis 
$\left\{\mathbf l,\mathbf e_1,\dots,\mathbf e_5\right\}$ of $\Pic S$).
Therefore $|D|$ contains a smooth connected member, that is a curve
$C \subset \mathbb P^4$ of degree $d$ and genus $g$
by \eqref{eqn:degree and genus}. We recall that all
smooth complete intersections $S_{2,2} \subset \mathbb P^4$ are
parametrised by an open subset $U$ of Grassmannian variety
$\Gr(2,V)$ of dimension $26$, 
where $V=H^0(\mathbb P^4,\mathcal O_{\mathbb P^4}(2))$.
Thus the pairs $(C,S)$ of $C$ and $S$ are parametrised by
$\mathcal U:=pr_2^{-1}(U)\cap \HF^{sc} \mathbb P^4$,
which is an open subset of $\mathbb P^{d+g-1}$-bundle over $U$, 
where $d+g-1=\dim |\mathcal O_S(D)|$.
Therefore curves $C \subset \mathbb P^4$ contained in 
some $S_{2,2} \subset \mathbb P^4$
with coordinate $(a;b_1,\dots,b_5)$ are parametrised by
$W=pr_1(\mathcal U)$, which is a 
locally closed irreducible subset of $\Hilb^{sc} \mathbb P^4$.
The relations between $U, \mathcal U$ and $W$ are illustrated by the diagram
$$
\begin{array}{ccccccc}
  \HF^{sc} \mathbb P^4 & \supset & \mathcal U & \mapright{pr_1} & W & \subset & \Hilb^{sc}  \mathbb P^4 \\
  && \mapdown{pr_2} &&&& \\
  \Gr(2,V) & \supset & U. &&&&
\end{array}
$$
It follows from the diagram that 
$\dim W=\dim \mathcal U=\dim U+(d+g-1)=d+g+25$.
\qed

\begin{rmk}
  \label{rmk:density and finiteness}
  By dimensions, the family $W$ in 
  Lemma~\ref{lem:maximal families and 6-tuples}
  is dense in the $S$-maximal family $W_{C,S}$
  for any member $C$ of $W$ and $S$ containing $C$.
  In fact, by the first projection $pr_1$, the unique irreducible component
  $\mathcal W_{C,S}$ of $\HF^{sc} \mathbb P^4$ passing through $(C,S)$ 
  is mapped onto the closure $\overline W$ of $W$ in $\Hilb^{sc} \mathbb P^4$.
  In particular, if we fix $d$ and $g$, then there exist only finitely many
  $S_{2,2}$-maximal families $\overline W \subset 
  \Hilb_{d,g}^{sc} \mathbb P^4$.
  One can also prove this fact by showing
  the number of $6$-tuples $(a;b_1,\dots,b_5)$ of integers
  satisfying \eqref{standard-prescribed-nef-big} is finite,
  together with Lemma~\ref{lem:maximal families and 6-tuples}
  (by using e.g.~Cauchy–Schwarz inequality on the coordinates $b_1,\dots,b_5$).
\end{rmk}

Here and later, we denote by $\Hilb_{d,g}^{sc} \mathbb P^4$
the subscheme of the Hilbert scheme $\Hilb^{sc} \mathbb P^4$
that parametrizes curves of degree $d$ and genus $g$.
Then according to Lemma~\ref{lem:maximal families and 6-tuples},
every $6$-tuple $(a;b_1,\dots,b_5)$ of integers 
satisfying \eqref{standard-prescribed-nef-big} corresponds to
a $S_{2,2}$-maximal family $W_{C,S_{2,2}}$ in $\Hilb_{d,g}^{sc} \mathbb P^4$.
In what follows, $W(a;b_1,\dots,b_5)$ represents
the corresponding $S_{2,2}$-maximal family in the Hilbert scheme.
By definition, $W(a;b_1,\dots,b_5)$ is irreducible and closed 
in $\Hilb_{d,g}^{sc} \mathbb P^4$ (and $\Hilb^{sc} \mathbb P^4$).

\begin{rmk}
  \label{rmk:dimension reason}
  By deformation theory, every irreducible component of 
  $\Hilb_{d,g}^{sc} \mathbb P^4$ is of dimension at least $5d+1-g$ 
  ($=\chi(C,N_{C/\mathbb P^4})$).
  Therefore, if $d >8$ and the $S_{2,2}$-maximal family $W_{C,S}$
  in $\Hilb^{sc} \mathbb P^4$
  is an irreducible component of $\Hilb_{d,g}^{sc} \mathbb P^4$,
  then $d+g+25 \ge 5d+1-g$, equivalently, $g \ge 2d-12$.
  Therefore, if $(d,g)$ belongs to 
  the interior of the region [III] or [IV] in Figure~\ref{fig:d-g-region},
  then $W_{C,S} \subset \Hilb_{d,g}^{sc} \mathbb P^4$ 
  is not a component, i.e.,~there exists an irreducible component $V$
  of $\Hilb^{sc} \mathbb P^4$ strictly containing $W_{C,S}$.
\end{rmk}

\paragraph{\bf Proof of Theorem~\ref{thm:main1}}
Let $C$ and $S$ be as in the theorem. Then $S$ is unique by $d >8$.
It follows from Lemma~\ref{lem:maximal families and 6-tuples} 
and Remark~\ref{rmk:density and finiteness}
that curves $C$ (of degree $d$ and genus $g$)
are parametrised by the union of open dense subsets $W$
of the $S_{2,2}$-maximal families
$W(a;b_1,\dots,b_5) \subset \Hilb^{sc} \mathbb P^4$ 
corresponding to $6$-tuples $(a;b_1,\dots,b_5)$ of integers
satisfying \eqref{standard-prescribed-nef-big}.
Moreover the number of such $6$-tuples $(a;b_1,\dots,b_5)$
is finite by the same remark, and hence we have proved (1).

For proving (2) and (3), we note by \eqref{isom:normal bundle of S} that
$N_{S/\mathbb P^4}(-C) \simeq \mathcal O_S(-D_2)^{\oplus 2}$,
where $D_2=C+2K_S$ is a divisor on $S$.
Thereby 
$H^i(N_{S/\mathbb P^4}(-C))
=H^i(-D_2)^{\oplus 2}$ for all $i$.
Then it follows from Lemma~\ref{lem:normality and speciality} that
$H^i(N_{S/\mathbb P^4}(-C))$ ($i=1,2$)
are isomorphic to $H^1(\mathcal I_C(2))^{\oplus 2}$
and $H^1(\mathcal O_C(2))^{\oplus 2}$, respectively.
If $C$ is $2$-normal, then $pr_1$ is smooth at $(C,S)$ 
(see ~\S\ref{subsec:flag and max_fam}),
and hence $\Hilb^{sc} \mathbb P^4$ is nonsingular at $[C]$
by Lemma~\ref{lem:nonsingular of expected dimension}.
This implies that $\Hilb^{sc} \mathbb P^4$ is smooth along $W$ and
its closure $\overline W$ ($=W(a;b_1,\dots,b_5)$)
is a generically smooth component of $\Hilb^{sc} \mathbb P^4$.
On the other hand, by \eqref{ineq:codimension of HF in Hilb}, 
there exists an inequality
$$
2h^1(\mathbb P^4,\mathcal I_C(2))-2h^1(C,\mathcal O_C(2))
\le 
\dim_{[C]} \Hilb^{sc} \mathbb P^4 - \dim_{(C,S)} \HF^{sc} \mathbb P^4
\le
2h^1(\mathbb P^4,\mathcal I_C(2)),
$$
where the inequality to the right is strict if and only if
$\Hilb^{sc} \mathbb P^4$ is singular at $[C]$.
Thus if $\mathcal O_C(2)$ is non-special, 
then $\Hilb^{sc} \mathbb P^4$ is still nonsingular along $W$
and $\overline W$ is of codimension $2h^1(\mathcal I_C(2))$
in $\Hilb^{sc} \mathbb P^4$.

Finally we prove (4). We see that if $g > (d^2-3d+6)/10$ 
then $C$ is $2$-normal by Lemma~\ref{lem:2-normal},
and $\mathcal O_C(2)$ is non-special
if either $2g-2=\deg K_C < \deg \mathcal O_C(2)=2d$, or
$g < 5d/2-35$ and $(d,g) \not\in P$ by Lemma~\ref{lem:2-nonspecial},
where $P$ is the set of exceptional pairs in \eqref{eqn:exceptional pairs}.
Thus we have completed the proof. \qed

\medskip

We remark that $S_{2,2}$-maximal families $W_{C,S}$ of curves $C \subset \mathbb P^4$
are not necessarily irreducible components of $\Hilb^{sc} \mathbb P^4$
even in the region [II] of Figure~\ref{fig:d-g-region}.
In fact, if $C$ is not linearly normal, then $W_{C,S}$ is not a component
by Example~\ref{ex:of type n} below.
Ellia~\cite{Ellia87} constructed a similar example
for $S_3$-maximal (or $3$-maximal) families of curves $C \subset \mathbb P^3$
contained in a smooth cubic surface, in order to
to give a counterexample to a conjecture of Kleppe (cf.~\cite[Conjecture 4]{Kleppe87})
and proposed a modification of his original conjecture
by putting an additional assumption of
the linear normality of $C \subset \mathbb P^3$.

\begin{ex}[cf.~\cite{Ellia87}]
  \label{ex:of type n}
  Let $S_5 \subset \mathbb P^5$ be a smooth quintic del Pezzo surface,~i.e.,
  the blow-up of $\mathbb P^2$ at general $4$ points.
  Let $\Gamma$ be a smooth complete intersection of $S_5$ 
  with a hypersurface of degree $n\ge 3$.
  Then $\Gamma$ has the coordinate $(3n;n,n,n,n)$ in $\Pic S_5 \simeq \mathbb Z^5$.
  Let $p_0 \in S_5\setminus \Gamma$ be a point not contained in any $10$ lines on $S_5$.
  Then the projection
  $\pi_{p_0}: \mathbb P^5 \dashrightarrow \mathbb P^4$ from the center $p_0$
  maps $\Gamma$ isomorphically onto a smooth curve $C_0$ in $\mathbb P^4$
  contained in a smooth complete intersection $S_{2,2}$ in $\mathbb P^4$,
  in which $C_0$ has the coordinate $(3n;n,n,n,n,0)$ in $\Pic S_{2,2} \simeq \mathbb Z^6$.
  Therefore there exists a diagram
  $$
  \begin{array}{ccccccc}
    \Gamma & \sim & (3n;n,n,n,n) & \subset & S_5 & \subset & \mathbb P^5 \\
    \mapdown{\simeq} &&&& \mapdown{} && \mapdown{\pi_{p_0}} \\
    C_0 & \sim & (3n;n,n,n,n,0) &\subset & S_{2,2} & \subset & \mathbb P^4. \\
  \end{array}
  $$
  We can check that $C_0$ is neither linearly normal nor quadratically normal
  (i.e.,~we have $h^1(\mathcal I_{C_0}(1))h^1(\mathcal I_{C_0}(2))\ne 0$).
  Then there exists a specialization $p_t \rightsquigarrow p_0$ of points in $\mathbb P^5$
  with $p_t \not\in S_5$ ($t \ne 0$), which induces specialization
  $C_t:=\pi_{p_t}(\Gamma) \rightsquigarrow C_0$ of curves $C_t$ in $\mathbb P^4$ to $C_0$.
  Here $C_t$ is isomorphic to $\Gamma$ and
  contained in a (non-normal) quintic surface $\pi_{p_t}(S_5)$ in $\mathbb P^4$.
  We note that $C_t$ is not contained in any $S_{2,2} \subset \mathbb P^4$,
  because otherwise, $C_t$ is contained in a complete intersection $Y_t$ 
  of $\pi_{p_t}(S_5)$ with a quadric hypersurface $Q \subset \mathbb P^4$
  by Lemma~\ref{lem:generic projection} below.
  However this is impossible by degree ($\deg(C_t)=5n > \deg(Y_t)=10$).
  By inverse construction, we see that all smooth curves 
  contained in $S_{2,2}$ with coordinate
  $(3n;n,n,n,n,0)$ are specializations of curves contained in 
  such a quintic surface in $\mathbb P^4$.
  Therefore, there exists an irreducible closed subset $V$
  of $\Hilb^{sc} \mathbb P^4$ that strictly contains $W_{C,S}$.
\end{ex}

\begin{lem}
  \label{lem:generic projection}
  The projection of every smooth quintic del Pezzo surface
  from a general point in $\mathbb P^5$ is
  not contained in any quadric hypersurface in $\mathbb P^4$.
\end{lem}
\Proof It is known that every smooth 
quintic del Pezzo surface $F \subset \mathbb P^5$
is defined by five $4\times 4$ Pfaffians of a $5\times 5$ 
skew-symmetric matrix of linear forms (cf.~e.g.~\cite{Buchsbaum-Eisenbud77}).
In particular, the homogeneous ideal $I_F$ of $F$ is generated
by $5$ quadratic forms $Q_i$ ($i=0,\dots,4$) on $\mathbb P^5$.
Thus the linear system $|\mathcal I_F(2)|$
of hyperquadrics containing $F$ is of dimension $4$.
Let $A_i$ ($i=0,\dots,4$) denote the $5\times 5$ symmetric matrices 
corresponding to $Q_i$, respectively. 
Then the locus of quadrics of rank at most $4$
(i.e.,~quadric cones) in $|\mathcal I_F(2)|\simeq \mathbb P^4$ 
is defined by the equation
\begin{equation}
  \label{eqn:locus of rank 4}
  \det (\lambda_0 A_0 +\dots+ \lambda_4 A_4)=0,
\end{equation}
where $\Lambda=(\lambda_0,\dots,\lambda_4)$ represents a point in $\mathbb P^4$.
Therefore, the general point $p$ of $\mathbb P^5$
is not contained in the set of 
all vertices of rank $4$ quadrics in $|\mathcal I_F(2)|$.
Similarly, the locus of rank (at most) $3$ quadrics in $|\mathcal I_F(2)|$
is defined by $5$ principal minors of $\lambda_0 A_0 +\dots+ \lambda_4 A_4$
together with \eqref{eqn:locus of rank 4} by the principal minor theorem.
Thus $p$ is not contained in the vertex lines of such quadrics. 
Since there exist no hyperquadrics of rank $2$ or $1$ in $|\mathcal I_F(2)|$,
there exist no hyperquadrics in $\mathbb P^4$ containing $\pi_p(F)$ 
if $p\in \mathbb P^5$ is general.
\qed 

\begin{lem}
  \label{lem:not 1-normal}
  Let $C$ be a smooth connected curve of degree $d$ in $\mathbb P^4$
  contained in a smooth complete intersection $S=S_{2,2}$ in $\mathbb P^4$.
  If $d>10$, $C$ is not linearly normal and 
  $C+K_S$ is effective and big,
  then there exists an irreducible component of $\Hilb^{sc} \mathbb P^4$
  which strictly contains $W_{C,S}$.
\end{lem}
\Proof
Since $H^1(\mathcal I_C(1))\simeq H^1(S,-C-K_S)$ by \eqref{isom:normality}
and $C+K_S$ is big, we see that
$C+K_S$ is not nef. This implies that there exists a line $E$ on $S$
such that $C.E=0$. Then by contraction,
there exists a smooth quintic del Pezzo surface $F \subset \mathbb P^5$
whose inner projection (i.e.,~a projection from a point on the surface)
is $S \subset \mathbb P^4$.
Furthermore, there exists a smooth curve $C'$ on $F$
mapped isomorphically onto $C$ by the projection.
Then by the same argument in Example~\ref{ex:of type n},
there exists a specialization $C_t \rightsquigarrow C$
of curves $C_t$ such that $C_0=C$ and
$C_t$ ($t\ne 0$) is not contained in any $S_{2,2} \subset \mathbb P^4$ 
by $\deg(C_t)=\deg(C)=d > 2 \cdot 5=2 \deg F$.
Thus $W_{C,S}$ is not a component of $\Hilb^{sc} \mathbb P^4$. 
\qed

\section{Obstructed curves and non-reduced components of $\Hilb^{sc} \mathbb P^4$}
\label{sec:non-reducedness}

In this section, we compute primary obstructions to deforming 
space curves lying on a smooth quartic del Pezzo surface and prove Theorem~\ref{thm:main2}.

\subsection{infinitesimal deformations with pole}
\label{subsec:with poles}

Let $S$ be a smooth del Pezzo surface and $E \subset S$ a line.
Then since $K_S.E=E^2=-1$, the invertible sheaf 
$\mathcal O_E(-K_S+E)$ on $E\simeq \mathbb P^1$ is trivial.
If $\deg S\ge 2$ then 
the complete linear system $|-K_S+E|$ on $S$ defines a morphism
$S \rightarrow \mathbb P^{\deg S+1}$,
which contract $E$ to a point on $S'$ (cf.~\cite[Theorem II.17]{Beauville})
and its image is again a smooth del Pezzo surface.
Here this morphism (or the image surface)
is called an {\em unprojection} of $S$.
Thus there exists a rational section of $\mathcal O_S(-K_S)$
admitting a pole only along $E$, and 
this section is unique up to scalar multiplications
together with modulo the global sections of $\mathcal O_S(-K_S)$.
We explicitly describe this section in the case $\deg S=4$
by using homogeneous coordinates of 
$\mathbb P^4$ and equations of $S$.

Let $S$ be a (smooth) complete intersection $S_{2,2} \subset \mathbb P^4$
and $E$ a line on $S$.
We take homogeneous coordinates $x_i$ ($0 \le i \le 4$) of $\mathbb P^4$,
so that the first three coordinates $x_0,x_1,x_2$ 
define $E$ in $\mathbb P^4$.
Then there exist two (linearly independent) quadratic polynomials 
\begin{equation}
  \label{eqn:S}
  q_i = x_0 a_{i0} +x_1 a_{i1} +x_2 a_{i2} \quad (i=1,2)
\end{equation}
defining $S$ in $\mathbb P^4$, where $a_{ij}$
are linear polynomials on $\mathbb P^4$ ($j=0,1,2$).
We consider a rational section $\theta$ of $\mathcal O_S(1)$
defined by
\begin{equation}
  \label{eqn:anticanonical section with pole}
  \theta = \dfrac{\Delta_0}{x_0} 
  = \dfrac{\Delta_1}{x_1} = \dfrac{\Delta_2}{x_2},
\end{equation}
where $\Delta_j$ ($j=0,1,2$) are $2$-minors of a $2\times 3$ matrix
\begin{equation}
  \label{eqn:matrix A}
  A=
  \begin{pmatrix}
    a_{10} & a_{11} & a_{12} \\
    a_{20} & a_{21} & a_{22}
  \end{pmatrix}.
\end{equation}
More precisely, $\Delta_j=(-1)^j \det A_j$, where
$A_j$ is the matrix obtained from $A$ by removing the column of index $j$.
We see that the six open subsets $D(x_j)$ and $D(\Delta_j)$ ($j=0,1,2$)
cover $S$ (and $\mathbb P^4$) by the smoothness of $S$,
where $D(f)$ denotes the standard (affine) open subset of $S$
associated to a section $f$ of an invertible sheaf on $S$
(on $\mathbb P^4$).
Here and later, we abuse notations and denote by the same symbols
$x_i$ (resp.~$\Delta_j$) the restrictions to $S$ of
the global sections of $\mathcal O_{\mathbb P^4}(1)$ 
(resp.~$\mathcal O_{\mathbb P^4}(2)$).
We denote the open covering 
$\left\{D(x_i),D(\Delta_i) \bigm| i=0,1,2 \right\}$ of $S$ 
(or $\mathbb P^4$) by $\mathfrak U_1$.

\begin{lem}
  The rational section $\theta$ of $\mathcal O_S(1)$ ($\sim -K_S$)
  (in \eqref{eqn:anticanonical section with pole})
  defines a global section of $\mathcal O_S(-K_S+E)$.
\end{lem}
\Proof
Since $x_0,x_1,x_2$ form a basis of 
$H^0(\mathbb P^4,\mathcal I_{E/\mathbb P^4}(1))$,
their restrictions to $S$ form a basis of $H^0(S,\mathcal O_S(1)(-E))$ ($\simeq k^3$).
Moreover, since $\mathcal O_S(1)(-E)\simeq \mathcal O_S(-K_S-E)$,
there exist global sections $x_0',x_1',x_2'$ of $\mathcal O_S(-K_S-E)$
corresponding to $x_0,x_1,x_2$, respectively.
Here we note that
\begin{equation}
  \label{eqn:x_i and x_i'}
  D(x_j)=D(x_j')\setminus E
\end{equation}
for all $j=0,1,2$.
Since the linear system $|-K_S-E|$ defines a morphism $S \rightarrow E$, sending
$p \mapsto (x_0'(p):x_1'(p):x_2'(p))$, corresponding to the projection
of $S$ from $E$, 
we see that $\mathfrak U_2:=\left\{D(x_0'),D(x_1'),D(x_2')\right\}$ covers whole $S$.
Thus $\theta$ defines a global section of $\mathcal O_S(-K_S+E)$
by its expression \eqref{eqn:anticanonical section with pole}.
\qed

\medskip

We remark that the correspondence 
\begin{equation}
  \label{isom:trivialization}
  \mathcal O_E \mapright{\sim} \mathcal O_E(-K_S+E), 
  \qquad
  \mu \longleftrightarrow \mu (\theta\big{\vert}_E)
\end{equation}
gives a trivialization. We define two $1$-cochains
$\hat t_i$ ($i=1,2$) of $\mathcal O_S(K_S+2E)$ by
\begin{equation}
  \label{eqn:fundamental 1-cocycles}
  \hat t_i:=\left(
    \dfrac{a_{i0}}{x_1x_2},
    \dfrac{a_{i1}}{x_2x_0},
    \dfrac{a_{i2}}{x_0x_1}
  \right)
  \in 
  C^1(\mathfrak U_2,\mathcal O_S(K_S+2E)).
\end{equation}
Here for each triple $(k,l,m)$ with $\left\{k,l,m\right\}=\left\{1,2,3\right\}$,
the symbol $a_{ik}/x_lx_m$ stands for the data of $\hat t_i$ over $D(x_l'x_m')$,
which is a local section of $\mathcal O_S(-1)(2E)\simeq \mathcal O_S(K_S+2E)$.
Then the $1$-cochains $\hat t_i$ ($i=1,2$) satisfy the cocycle conditions,
i.e.,~$(\hat t_i)_{12}+(\hat t_i)_{20}+(\hat t_i)_{01}=0$ by \eqref{eqn:S}.
Thus $\hat t_i$ represents an element $t_i$ in $H^1(S,\mathcal O_S(K_S+2E))$ for $i=1,2$.
In what follows, by abuse of notations,
we denote a cohomology class and a cochain representing it
by the same symbol. 
Thus we write 
$t_i:=\left(a_{i0}/x_1x_2,a_{i1}/x_2x_0,a_{i2}/x_0x_1\right)
\in H^1(S,\mathcal O_S(K_S+2E))$.
We see that the restriction of $\mathcal O_S(K_S+2E)$ to $E$
is isomorphic to $\mathcal O_E(-3)$.
The following lemma will be a key to our proof of Theorem~\ref{thm:obstruction} later.

\begin{lem}
  \label{lem:linear independence of coh classes}
  Under the trivialization \eqref{isom:trivialization}, we have
  $$
  t_i\big{\vert}_E
  =\left(
    \dfrac{a_{i0}}{\Delta_1\Delta_2},
    \dfrac{a_{i1}}{\Delta_2\Delta_0},
    \dfrac{a_{i2}}{\Delta_0\Delta_1}
  \right)
  \in 
  H^1(\mathfrak U_2,\mathcal O_E(-3))
  $$
  for $i=1,2$. Moreover, $t_1\big{\vert}_E$ and $t_2\big{\vert}_E$
  are linearly independent, and hence they form a basis of
  $H^1(E,\mathcal O_E(-3))\simeq k^2$.
\end{lem}
\Proof
There exists a natural exact sequence 
\begin{equation}
  \label{ses:normal bundles}
  \begin{CD}
    0 @>>> \underbrace{N_{E/S}}_{\simeq \mathcal O_E(-1)} 
    @>{\iota}>> 
    \underbrace{N_{E/\mathbb P^4}}_{\simeq \mathcal O_E(1)^{\oplus 3}}
    @>{\pi_{E/S}}>> 
    \underbrace{N_{S/\mathbb P^4}\big{\vert}_E}_{\simeq \mathcal O_E(2)^{\oplus 2}}
    @>>> 0
  \end{CD}
\end{equation}
on $E \simeq \mathbb P^1$. Taking the dual bases of 
$x_0,x_1,x_2$ and $q_1,q_2$, we see that
the sheaf homomorphisms $\pi_{E/S}$ and $\iota$
are represented by the matrix $A$ in \eqref{eqn:matrix A} and 
${}^t \begin{pmatrix}
  \Delta_0 & \Delta_1 & \Delta_2
\end{pmatrix}
$, respectively. Here and later, ${}^t M$ denotes the transpose of a matrix $M$.
The exact sequence obtained by tensoring $\eqref{ses:normal bundles}$ 
with $\mathcal O_E(-2)$ induces a coboundary map
$$
\partial_{E/S}: H^0(E,
\underbrace{N_{S/\mathbb P^4}\big{\vert}_E\otimes \mathcal O_E(-2)}_{\simeq \mathcal O_E^{\oplus 2}}) \rightarrow H^1(E,\underbrace{N_{E/S}\otimes \mathcal O_E(-2)}_{\mathcal O_E(-3)}),
$$
which is clearly an isomorphism.
We will show that the global sections
$e_1:=(1, 0)$ and $e_2:=(0, 1)$ of the trivial bundle $\mathcal O_E^{\oplus 2}$
are respectively mapped to 
$-t_2 \big{\vert}_E$ and $t_1 \big{\vert}_E$ by $\partial_{E/S}$.

First we consider $e_1 \in H^0(\mathcal O_E^{\oplus 2})$. 
We define a local section $s_i$ ($i=0,1,2$) of $\mathcal O_E(-1)^{\oplus 3}$ 
over $D(x_i')$ by
$$
s_0: = \dfrac{1}{\Delta_0}
\begin{pmatrix}
  0 \\ a_{22} \\ -a_{21}
\end{pmatrix},
s_1: = \dfrac{1}{\Delta_1}
\begin{pmatrix}
  -a_{22} \\  0\\ a_{20}
\end{pmatrix},
s_2: = \dfrac{1}{\Delta_2}
\begin{pmatrix}
  a_{21} \\  -a_{20} \\ 0
\end{pmatrix}.
$$
Then by computations, we can show that 
$\pi_{E/S}(s_i)=e_1$ for all $i=0,1,2$. For example, we see that
$$
\pi_{E/S}(s_0)
=\dfrac{1}{\Delta_0}
A\begin{pmatrix}
  0 \\ a_{22} \\ -a_{21}
\end{pmatrix}
=\dfrac{1}{\Delta_0}
\begin{pmatrix}
  a_{11}a_{22}-a_{12}a_{21} \\ 0
\end{pmatrix}
=\begin{pmatrix}
  1 \\ 0
 \end{pmatrix}
=e_1.
$$
This implies that $s_i$ ($i=0,1,2$) are local lifts of $e_1$ by $\pi_{E/S}$
over $D(x_i')$.
Then since $a_{20}\Delta_0+a_{21}\Delta_1+a_{22}\Delta_2=0$, we have
\begin{align*}
  \iota (-t_2\big{\vert}_E)
  & =\left(
      \dfrac{-a_{20}}{\Delta_1\Delta_2},
      \dfrac{-a_{21}}{\Delta_2\Delta_0},
      \dfrac{-a_{22}}{\Delta_0\Delta_1}
    \right)
    \begin{pmatrix}
      \Delta_0 \\ \Delta_1 \\ \Delta_2
    \end{pmatrix}
  =\begin{pmatrix}
    \left(
      \dfrac{-a_{20}\Delta_0}{\Delta_1\Delta_2},
      \dfrac{-a_{21}}{\Delta_2},
      \dfrac{-a_{22}}{\Delta_1}
    \right)
    \\
    \left(
      \dfrac{-a_{20}}{\Delta_2},
      \dfrac{-a_{21}\Delta_1}{\Delta_2\Delta_0},
      \dfrac{-a_{22}}{\Delta_0}
    \right)
    \\
    \left(
      \dfrac{-a_{20}}{\Delta_1},
      \dfrac{-a_{21}}{\Delta_0},
      \dfrac{-a_{22}\Delta_2}{\Delta_0\Delta_1}
    \right)
  \end{pmatrix}\\
  &=
  \begin{pmatrix}
    \left(
      \dfrac{a_{22}}{\Delta_1}+\dfrac{a_{21}}{\Delta_2},
      \dfrac{-a_{21}}{\Delta_2},
      \dfrac{-a_{22}}{\Delta_1}
    \right)
    \\
    \left(
      \dfrac{-a_{20}}{\Delta_2},
      \dfrac{a_{22}}{\Delta_0}+\dfrac{a_{20}}{\Delta_2},
      \dfrac{-a_{22}}{\Delta_0}
    \right)
    \\
    \left(
      \dfrac{-a_{20}}{\Delta_1},
      \dfrac{-a_{21}}{\Delta_0},
      \dfrac{a_{21}}{\Delta_0}+\dfrac{a_{20}}{\Delta_1}
    \right)
  \end{pmatrix}
  =(s_2-s_1,s_0-s_2,s_1-s_0),
\end{align*}
and hence $\iota(-t_2\big{\vert}_E)$ is equal to 
the coboundary of the $0$-cochain $(s_0,s_1,s_2)$ 
of $\mathcal O_E(-1)^{\oplus 3}$.
This implies that $\partial_{E/S}(e_1)=-t_2\big{\vert}_E$.
Similarly, we can prove $\partial_{E/S}(e_2)=t_1\big{\vert}_E$
by computations. Since $e_1$ and $e_2$ are linearly independent,
so are $t_1\big{\vert}_E$ and $t_2\big{\vert}_E$.
\qed

We prepare a lemma concerning infinitesimal deformations with pole.

\begin{lem}
  \label{lem:2-form on S}
  Let $\beta$ be a global section of $N_{S/\mathbb P^4}(E)$,~i.e.,
  an infinitesimal deformation of $S$ with pole along $E$.
  Then there exist a global section $\xi$ of 
  $\mathcal O_{\mathbb P^4}(1)^{\oplus 2}$
  and a global section $\eta$ of $\mathcal O_{\mathbb P^4}(2)^{\oplus 2}$
  such that
  $$\beta=\theta(\xi\big{\vert}_S)+\eta\big{\vert}_S.$$
\end{lem}
\Proof
There exists an exact sequence 
$$
\left[0 \longrightarrow \mathcal O_S(-E)
  \mapright{\scriptsize
    \begin{pmatrix}
     -1 & \theta 
    \end{pmatrix}}
    \mathcal O_S\oplus \mathcal O_S(-K_S)
    \mapright{\scriptsize
  \begin{pmatrix}
    \theta \\ 1
  \end{pmatrix}}
  \mathcal O_S(-K_S+E)
  \longrightarrow 0\right]\otimes \mathcal O_S(-K_S)
$$
on $S$ of Koszul type. Then since $H^1(-K_S-E)=0$, 
the induced map on cohomology groups from $H^0(-K_S)\oplus H^0(-2K_S)$
to $H^0(-2K_S+E)$ by ${}^t (\theta \ 1)$ is surjective.
Then the lemma follows from that
$N_{S/\mathbb P^4}(E)\simeq \mathcal O_S(-2K_S+E)^{\oplus 2}$
and the projectively normality of $S$.
In fact, the restriction map
$H^0(\mathcal O_{\mathbb P^4}(n)) \rightarrow 
H^0(\mathcal O_S(n)) \simeq H^0(-n K_S)$ is surjective for all 
$n \in \mathbb Z$.
\qed

\subsection{Obstructions}
\label{subsec:obstructions}

We devote the whole section to the proof of the following theorem,
a refined version of Theorem~\ref{thm:main2}.

\begin{thm}
  \label{thm:obstruction}
  Let $S$ be a smooth complete intersection $S_{2,2} \subset \mathbb P^4$,
  $C$ a smooth connected curve on $S$ such that
  \begin{enumerate}
    \renewcommand{\labelenumi}{{\rm (\roman{enumi})}}
    \item $C+3K_S \ge 0$, and
    \item the fixed part of $|C+2K_S|$ is a (single) line $E$ on $S$.
  \end{enumerate}
  Then 
  \begin{enumerate}
    \item the cokernel of the tangent map
    $$
    p_1: T_{\HF \mathbb P^4,(C,S)} 
    \rightarrow T_{\Hilb \mathbb P^4,[C]}=H^0(C,N_{C/\mathbb P^4})
    $$
    of the first projection $pr_1$ at $(C,S)$ is of dimension $2$.
    \item If $C$ is general in the class $[C]$ of $\Pic S$, then
    $\ob(\alpha)\ne 0$ for every 
    $\alpha \in H^0(N_{C/\mathbb P^4})\setminus \im p_1$.
  \end{enumerate}
\end{thm}

Theorem~\ref{thm:obstruction} says that there exist first order
infinitesimal deformations $\tilde C$ of $C$ in $\mathbb P^4$
not contained in any first order deformations $\tilde S$ of $S$ in 
$\mathbb P^4$. However, they do not lift to 
deformations ${\tilde {\tilde C}}$ of $C$ over $k[t]/(t^3)$,
and hence nor to any global deformations of $C$ in $\mathbb P^4$.

Let $C,S$ be as in Theorem~\ref{thm:obstruction}.
In what follows, for the sake of conventions,
we denote by $L$ the class of an invertible sheaf $\mathcal O_S(C+2K_S)$ 
on $S$,~i.e.,~$L=[\mathcal O_S(C+2K_S)]$ in $\Pic S$, 
and for a divisor $D$ on $S$, we denote the class of 
$L \otimes_S \mathcal O_S(D)$ by $L+D$, similarly.
Then $L+K_S\ge 0$ by (i) and hence $L \ge 0$.
It follows from (ii) that
$$
L=|L-E|+E,
$$
where $L-E$ is nef. This implies that $L.E=-1$ and
$E$ is the only line on $S$ with $L.E<0$.
Then we see $C.E=(L-2K_S).E=1$ and 
the scheme-theoretic intersection $C \cap E$ is just a point.

\begin{claim}
  \label{claim:positivity and vanishing}
  $\chi(S,-L)\ge 0$, $L^2>0$, 
  and $H^i(S,-L+E)=0$ for $i=0,1$.
  Moreover, if $C$ is general in $|C|$, then $p=C\cap E$ 
  is a general point on $E\simeq \mathbb P^1$.
\end{claim}
\Proof
Since $L-E$ is nef and $L+K_S\ge 0$,
we have by the Riemann-Roch theorem on $S$ that
$$
2\chi(S,-L)-2=L.(L+K_S)\ge E.(L+K_S)=-2
$$
and hence $\chi(S,-L)\ge 0$.
Then $L$ is big by \S\ref{subsec:delpezzo}({\bf P6})
and so is $L-E$ by $(L-E)^2=L^2+1$.
Then Kawamata-Viehweg vanishing theorem shows
$H^i(S,-L+E)=0$ for $i=0,1$.
Finally, the last assertion follows from $H^1(S,C-E)=0$.
Indeed, this implies
the restriction map $|C| \dashrightarrow |\mathcal O_E(1)|$
sending $C' \in |C|$ to $C'\cap E$ is dominant.
\qed

We are ready to prove the first statement of Theorem~\ref{thm:obstruction}.

\medskip

\paragraph{\bf Proof of Theorem~\ref{thm:obstruction}~(1)}

It follows from Lemma~\ref{lem:nonsingular of expected dimension}
and \eqref{ses:first tangent map} that
there exists an exact sequence
\begin{equation}
  \label{ses:cokernel of p_1}
  H^0(\mathbb P^4,N_{(C,S)/\mathbb P^4})
  \overset{p_1}{\longrightarrow} H^0(C,N_{C/\mathbb P^4})
  \overset{\delta}{\longrightarrow}
  H^1(S,N_{S/\mathbb P^4}(-C))
  \longrightarrow 0.
\end{equation}
Since $N_{S/\mathbb P^4}(-C)\simeq -L^{\oplus 2}$
by \eqref{isom:normal bundle of S}, we have
$H^i(S,N_{S/\mathbb P^4}(E-C))=0$ for $i=0,1$
by Claim~\ref{claim:positivity and vanishing}.
Moreover, it follows from
$H^1(S,N_{S/\mathbb P^4})=H^1(S,-2K_S^{\oplus 2})=0$
that the restriction map $H^0(S,N_{S/\mathbb P^4}(E))
\rightarrow H^0(E,N_{S/\mathbb P^4}(E)\big{\vert}_E)$ is surjective.
Then applying the snake lemma to
the commutative diagram \eqref{diag:comm1} with 
$\mathcal F=N_{S/\mathbb P^4}$,
we see that $H^1(S,N_{S/\mathbb P^4}(-C))$ is isomorphic 
to $H^0(E,N_{S/\mathbb P^4}(E-C)\big{\vert}_E)
\simeq H^0(E,\mathcal O_E(-L+E)^{\oplus 2})$,
where $\mathcal O_E(-L+E)$ is trivial.
Thus we see that $\coker p_1$ is of dimension $2$. \qed

\medskip

As we saw in \S\ref{subsec:primary and exterior},
every element $\alpha$ in $H^0(N_{C/\mathbb P^4})$
($\simeq \Hom_{C}(N_{C/\mathbb P^4}^{\vee},\mathcal O_C)
\simeq \Hom_{\mathbb P^4}(\mathcal I_C,\mathcal O_C)$)
defines a global $\mathcal O_C$-linear homomorphism
$\alpha': N_{C/\mathbb P^4}^{\vee} \rightarrow \mathcal O_C$ 
and $\ob(\alpha)$ is given by 
$\ob(\alpha)=\alpha' \cup \mathbf e' \cup \alpha'$,
where $\mathbf e' \in
\Ext^1_{\mathbb P^4}(\mathcal O_C,N_{C/\mathbb P^4}^{\vee})$
is the extension class of the exact sequence
\begin{equation}
  \label{ses:conormal2}
  0 \rightarrow N_{C/\mathbb P^4}^{\vee} \rightarrow 
  \mathcal O_{\mathbb P^4}/\mathcal I_C^2
  \rightarrow \mathcal O_C \rightarrow 0
\end{equation}
on $\mathbb P^4$.
Then by functoriality, for any integer $l$, $\ob(\alpha)$ 
induces a map
$$
\cup \ob(\alpha): 
H^0(C,N_{C/\mathbb P^4}^{\vee}(l)) \rightarrow H^1(C,\mathcal O_C(l))
$$
on cohomology groups on $C$.
Then by definition, the map has a decomposition
$$
  \begin{CD}
    H^0(N_{C/\mathbb P^4}^{\vee}(l))
    @>{\cup \alpha'}>>
    H^0(\mathcal O_C(l))
    @>{\cup \mathbf e'}>>
    H^1(N_{C/\mathbb P^4}^{\vee}(l))
    @>{\cup \mathbf \alpha'}>>
    H^1(\mathcal O_C(l))
  \end{CD}
$$
into three cup product maps on $\mathbb P^4$ (or $C$),
where $\cup \mathbf e'$ is the first coboundary map of
$\eqref{ses:conormal2}\otimes \mathcal O_{\mathbb P^4}(l)$.
We note that $\cup \alpha'$ and $\cup \mathbf e'$ preserve
the degree $l$, and so does their composition $\ob(\alpha)$.
We observe $N_{S/\mathbb P^4}\big{\vert}_C \simeq \mathcal O_C(2)^{\oplus 2}$.
Then by taking the direct sum 
of two copies of the above sequence with $l=2$,
we define a map $\Phi(\alpha)$ from $H^0(N_{C/\mathbb P^4}^{\vee}(2)^{\oplus 2})$
to $H^1(N_{S/\mathbb P^4}\big{\vert}_C)$ by the compositions
$$
  \begin{CD}
    \Phi(\alpha):
    H^0(N_{C/\mathbb P^4}^{\vee}(2)^{\oplus 2})
    @>{\cup \alpha'}>>
    H^0(N_{S/\mathbb P^4}\big{\vert}_C)
    @>{\cup \mathbf e'}>>
    H^1(N_{C/\mathbb P^4}^{\vee}(2)^{\oplus 2})
    @>{\cup \mathbf \alpha'}>>
    H^1(N_{S/\mathbb P^4}\big{\vert}_C).
  \end{CD}
$$
Then, in order to prove $\ob(\alpha)\ne 0$ for
$\alpha \in H^0(N_{C/\mathbb P^4})\setminus \im p_1$, 
it suffices to show $\Phi(\alpha)\ne 0$.
Taking the dual of 
$\pi_{C/S}:N_{C/\mathbb P^4} \twoheadrightarrow 
N_{S/\mathbb P^4}\big{\vert}_C$
and tensoring $\mathcal O_C(2)$, we obtain an inclusion
\begin{equation}
  \label{inclusion:natural}
  N_{S/\mathbb P^4}\big{\vert}_C^{\vee}(2) \simeq \mathcal O_C^{\oplus 2}
  \hookrightarrow N_{C/\mathbb P^4}^{\vee}(2).
\end{equation}
of sheaves on $C$. Thus we observe that $N_{C/\mathbb P^4}^{\vee}(2)$
contains a trivial subbundle $\mathcal O_C^{\oplus 2}$ on $C$.
We use this fact for the computation of $\Phi(\alpha)$ later.

Let $q_1,q_2$ be the two quadratic polynomials defining $S$ in $\mathbb P^4$,
which we have defined in \S\ref{subsec:with poles}.
Then their images in $H^0(N_{C/\mathbb P^4}^{\vee}(2))$, 
denoted by $\bar q_1,\bar q_2$, form a free basis of
the subbundle
$N_{S/\mathbb P^4}\big{\vert}_C^{\vee}(2) \simeq \mathcal O_C^{\oplus 2}$.
Let $\alpha \in H^0(N_{C/\mathbb P^4})\setminus \im p_1$
be an arbitrary element, and put
$u_i:=\alpha(q_i)$ in $H^0(\mathcal O_C(2))$ ($i=1,2$).
Since $\mathcal I_S$ is globally generated by $q_i$ ($i=1,2$),
the restriction of $\alpha: \mathcal I_C \rightarrow \mathcal O_C$
to $\mathcal I_S \subset \mathcal I_C$ is determined by $u_i$ ($i=1,2$).
Thus we identify
$\pi_{C/S}(\alpha) \in H^0(N_{S/\mathbb P^4}\big{\vert}_C)$
with $(u_1,u_2) \in H^0(\mathcal O_C(2)^{\oplus 2})$.
Then since $\alpha(q_i)=\alpha'(\bar q_i)$ ($i=1,2$), we have
\begin{equation}
  \label{eqn:first product}
  (\bar q_1,\bar q_2) \cup \alpha'
  =(\alpha(q_1),\alpha(q_2))=(u_1,u_2)=\pi_{C/S}(\alpha).
\end{equation}
Similarly, we identify the exterior component 
$\ob_S(\alpha)$ 
in $H^1(N_{S/\mathbb P^4}\big{\vert}_C)$
with $\Phi(\alpha)(\bar q_1,\bar q_2)$
in $H^1(\mathcal O_C(2)^{\oplus 2})$ by
\begin{equation}
  \label{eqn:second product}
  \ob_S(\alpha)=\pi_{C/S}(
    \alpha' \cup \mathbf e' \cup \alpha')
  =\pi_{C/S}
  (\alpha) \cup \mathbf e' \cup \alpha'
  =\Phi(\alpha)(\bar q_1,\bar q_2).
\end{equation}
It follows from the commutative diagram
$$
\begin{CD}
  N_{(C,S)/\mathbb P^4} @>{\pi_2}>> N_{S/\mathbb P^4}  \\
  @V{\pi_1}VV @V{|_C}VV \\
  N_{C/\mathbb P^4} @>{\pi_{C/S}}>> N_{S/\mathbb P^4}\big{\vert}_C
\end{CD}
$$
with \eqref{ses:normal1} and \eqref{ses:normal2}, or more directly,
\cite[Lemma~2.2]{Nasu6} that
the coboundary map $\delta$ in \eqref{ses:cokernel of p_1}
factors through the coboundary map
$$
\cup \mathbf k_C: 
H^0(C,N_{S/\mathbb P^4}\big{\vert}_C) \rightarrow H^1(S,N_{S/\mathbb P^4}(-C))$$
(see \eqref{ses:CE} for $\mathbf k_C$).
Then by assumption, we have
$\pi_{C/S}(\alpha)\cup \mathbf k_C\ne 0$.
We note that $H^1(S,N_{S/\mathbb P^4}(E-C))=0$ by Claim~\ref{claim:positivity and vanishing}.
Then Lemma~\ref{lem:restriction to C and E} shows that
$r(\pi_{C/S}(\alpha),E)$ lifts to a global section $\beta$ of 
$N_{S/\mathbb P^4}(E)$, but not that of $N_{S/\mathbb P^4}$,
i.e.,~an infinitesimal deformation of $S \subset \mathbb P^4$ with pole along $E$.
Moreover, since $\pi_{C/S}(\alpha) \cup \mathbf k_C \ne 0$,
the same lemma shows that the principal part
\begin{equation}
  \label{eqn:principal part}
  \beta\big{\vert}_E
  \in H^0(E,N_{S/\mathbb P^4}(E)\big{\vert}_E)
\end{equation}
of $\beta$ along $E$ is nonzero and it comes from a global section of
$N_{S/\mathbb P^4}(E-C)\big{\vert}_E$.

The following lemma is another key to our proof of Theorem~\ref{thm:obstruction}~(2).
In this lemma, we fix the equations
$x_0,x_1,x_2$ defining $E$ 
together with $q_1,q_2$ defining $S$ as in \S\ref{subsec:with poles}.
We consider the {\v C}ech cohomology groups 
$H^i(\mathfrak U_2,*)$ on $C$ (or $S$)
with respect to the open affine covering
$\mathfrak U_2:=\left\{D(x_j') \bigm| j=0,1,2 \right\}$ of $C$ (or $S$),
where $x_j'$ are global sections of $\mathcal O_S(-K_S-E)$
corresponding to $x_j$ ($j=0,1,2$) (see \S\ref{subsec:with poles}).
Let $t_i$ ($i=1,2$) be the elements of 
$H^1(S,\mathcal O_S(-1)(2E))$
represented by \eqref{eqn:fundamental 1-cocycles}.
Then we define an element $\mathbf t$ of $H^1(S,\mathcal O_S(-1)(2E)^{\oplus 2})$
by $\mathbf t:=(-t_2,t_1)$.
\begin{lem}
  \label{lem:key lemma}
  Let $C,S,E,p_1,\alpha$ be as in Theorem~\ref{thm:obstruction},
  and put $p:=C\cap E$.
  Let $\pi_{C/S}(\alpha)$ and $\ob_S(\alpha)$
  be the exterior components of $\alpha$ and $\ob(\alpha)$, respectively
  (cf.~\S\ref{subsec:primary and exterior}).
  \begin{enumerate}
    \item Put $\pi_{C/S}(\alpha):=(u_1,u_2)$, where 
    $u_i$ ($i=1,2$) are the $i$-th components of 
    $N_{S/\mathbb P^4}\big{\vert}_C \simeq \mathcal O_C(2)^{\oplus 2}$.
    Then for every $i=1,2$, the reduction $r(u_i \cup \mathbf e',2p)$
    is contained in a smaller subgroup
    $$
    H^1(C,\mathcal O_C(2p)^{\oplus 2}) 
    \subset H^1(C,N_{C/\mathbb P^4}^{\vee}(2)(2p)),
    $$
    where $\mathcal O_C^{\oplus 2}(2p)$
    is regarded as a subsheaf of
    $N_{C/\mathbb P^4}^{\vee}(2)(2p)$ by 
    $\eqref{inclusion:natural}\otimes \mathcal O_C(2p)$.
    \item 
    The reduction $r(\ob_S(\alpha),2p)$ of $\ob_S(\alpha)$
    in $H^1(N_{S/\mathbb P^4}\big{\vert}_C(2p))$
    is equal to the cup product
    \begin{equation}
      \label{eqn:exterior component}
      \xi \left(
      \mathbf t\big{\vert}_C \cup \pi_{C/S}(\alpha)
      \right),
    \end{equation}
    where 
    \begin{itemize}
      \item $\mathbf t\big{\vert}_C=(-t_2\big{\vert}_C,t_1\big{\vert}_C)$
      in $H^1(\mathcal O_C(-1)(2p)^{\oplus 2})
      \simeq H^1(N_{S/\mathbb P^4}\big{\vert}_C^{\vee}(1)(2p))$
      \item $\cup$ is taken by the map
      $$
      H^1(C,N_{S/\mathbb P^4}\big{\vert}_C^{\vee}(1)(2p))\times
      H^0(C,N_{S/\mathbb P^4}\big{\vert}_C)
      \overset{\cup}{\longrightarrow} 
      H^1(C,\mathcal O_C(1)(2p))
      $$
      \item $\xi$ is a global section of 
      $\mathcal O_{\mathbb P^4}^{\oplus 2}(1)$
      such that $\xi\big{\vert}_E$ is nonzero and has a zero at $p$.
    \end{itemize}
  \end{enumerate}
\end{lem}
\Proof The proof is similar to that of \cite[Lemma~3.5]{Nasu1}.
Its strategy is as follows.
We will compute the image of $u_i$ ($i=1,2$)
by the coboundary map
$$
H^0(C,\mathcal O_C(2)) \overset{\cup \mathbf e'}{\longrightarrow} 
H^1(C,N_{C/\mathbb P^4}^{\vee}(2))
$$
of $\eqref{ses:conormal2}\otimes_{\mathbb P^4} \mathcal O(2)$, using
the {\v C}ech cohomology with respect to the open affine covering
$\mathfrak U_1=\left\{D(x_j),D(\Delta_j) \bigm| j=0,1,2 \right\}$
of $\mathbb P^4$ (cf.~\S\ref{subsec:with poles}).
Admitting a pole of order $2$ at $p$,
we will see that a $1$-cocycle of 
$N_{S/\mathbb P^4}\big{\vert}_C^{\vee}(2)\simeq \mathcal O_C^{\oplus 2}$
represents the coboundary image $u_i\cup \mathbf e'$
 (cf.~\eqref{inclusion:natural}).
By considering another open affine covering $\mathfrak U_3$ of $C$ and
using isomorphisms between {\v C}ech cohomology groups
$H^1(\mathfrak U_k,\mathcal O_C^{\oplus 2}(2p))$
($k=1,2,3$), we will further show that $r(u_i\cup \mathbf e', 2p)$ is contained in
$H^1(\mathfrak U_2,\mathcal O_C^{\oplus 2}(2p))$.
Finally we will relate the cohomology classes
$r(u_i\cup \mathbf e', 2p)$ ($i=1,2$) and $r(\ob_S(\alpha),2p)$
to the cup product \eqref{eqn:exterior component} in the statement.

We first recall that $r(\pi_{C/S}(\alpha),p)$ lifts to
the global section $\beta$ of $N_{S/\mathbb P^4}(E)$
and its principal part $\beta\big{\vert}_E$, that is a global section of
$N_{S/\mathbb P^4}(E)\big{\vert}_E\simeq \mathcal O_E(1)^{\oplus 2}$,
is nonzero and has a zero at $p$ (cf.~\eqref{eqn:principal part}).
Then Lemma~\ref{lem:2-form on S} shows that
there exist a global sections 
$\xi=(\xi_1,\xi_2)$ of $\mathcal O_{\mathbb P^4}(1)^{\oplus 2}$ and 
$\eta =(\eta_1,\eta_2)$ of $\mathcal O_{\mathbb P^4}(2)^{\oplus 2}$
such that $\beta=\xi \theta + \eta$,
where $\theta$ is a global section of the sheaf $\mathcal O_S(1)(E)$ on $S$
defined by \eqref{eqn:anticanonical section with pole}.
Then via the trivialization \eqref{isom:trivialization},
we obtain $\beta\big{\vert}_E=\xi\big{\vert}_E$ as a section of
$N_{S/\mathbb P^4}(E)\big{\vert}_E\simeq \mathcal O_E(1)^{\oplus 2}$, 
which is nonzero and has a zero at $p$. Moreover, we have
$$
r(\pi_{C/S}(\alpha),p)=\beta\big{\vert}_C
=\xi\big{\vert}_C(\theta\big{\vert}_C)+\eta\big{\vert}_C,
$$
where $\theta\big{\vert}_C$ is a global section of $\mathcal O_C(1)(p)$.
Therefore, in terms of the linear forms $\xi_i$ and
quadratic forms $\eta_i$ on $\mathbb P^4$, the global sections $u_i$ ($i=1,2$)
of $\mathcal O_C(2)$ are expressed as
\begin{equation}
  \label{eqn:u_ij}
  u_{ij}= \dfrac{\xi_i \Delta_j}{x_j} + \eta_i
\end{equation}
over the open set $D(x_j)$ ($j=0,1,2$) of $C$ (or $\mathbb P^4$).
Then the coboundary image 
$u_i \cup \mathbf e'$ of $u_i$ in $H^1(\mathfrak U_1,N_{C/\mathbb P^4}^{\vee}(2))$
is represented by
$$
u_{i1}-u_{i0}=\xi_i\left(\dfrac{ \Delta_1}{x_1}-\dfrac{\Delta_0}{x_0}\right)
$$
over $D(x_0)\cap D(x_1)$ and we compute that
\begin{align*}
  \dfrac{\Delta_1}{x_1}-\dfrac{\Delta_0}{x_0}
  &=\dfrac{x_0(a_{12}a_{20}-a_{10}a_{22})
    -x_1(a_{11}a_{22}-a_{12}a_{21})}{x_0x_1}\\
  &=\dfrac{-a_{22}(x_0a_{10}+x_1a_{11})+a_{12}(x_0a_{20}+x_1a_{21})}{x_0x_1}\\
  &=\dfrac{-a_{22}\bar q_1+a_{12}\bar q_2}{x_0x_1},
\end{align*}
where the last equality follows from \eqref{eqn:S}.
Similarly, it follows from \eqref{eqn:u_ij} that
$u_i \cup \mathbf e'$ is represented by
$$
\xi_i\left(
  \dfrac{\Delta_2}{x_2}-\dfrac{\Delta_1}{x_1}
\right)=
\xi_i\left(
  \dfrac{-a_{20}\bar q_1+a_{10}\bar q_2}{x_1x_2}
\right)
\quad
\mbox{and}
\quad
\xi_i\left(
  \dfrac{\Delta_0}{x_0}-\dfrac{\Delta_2}{x_2}
\right)=
\xi_i\left(
  \dfrac{-a_{21}\bar q_1+a_{11}\bar q_2}{x_2x_0}
\right)
$$
over $D(x_1)\cap D(x_2)$ and $D(x_0)\cap D(x_2)$, respectively.
The above computations shows
that over the subsets $\mathfrak U_1':=\left\{ D(x_j) \bigm| j=0,1,2 \right\}$
of $\mathfrak U_1$, $u_i \cup \mathbf e'$ is respectively represented by
$$
\xi_i\left(
-\hat t_2\big{\vert}_C \bar q_1 + \hat t_1\big{\vert}_C \bar q_2
\right),
$$
where $\hat t_i$ ($i=1,2$) are the $1$-cocycles of $\mathcal O_S(-1)(2E)$
representing $t_i$ (cf.~\eqref{eqn:fundamental 1-cocycles}).
This implies that with respect to the subcovering $\mathfrak U_1'$
of $\mathfrak U_1$, 
$u_i \cup \mathbf e'$ is represented by a $1$-cocycle of
$
N_{S/\mathbb P^4}\big{\vert}_C^{\vee}(2)
\subset N_{C/\mathbb P^4}^{\vee}(2)
$.

On the other hand, we recall by \eqref{eqn:x_i and x_i'} that
the subcovering $\mathfrak U_1' \subset \mathfrak U_1$
covers whole $S$ except for $E$.
Then $\mathfrak U_1'$ covers whole $C$ except for $p=C \cap E$ and
we have $D(x_j)=D(x_j')\setminus p$ ($j=0,1,2$) over $C$.
Therefore, $\xi_i (- \hat t_2\big{\vert}_C, \hat t_1\big{\vert}_C)$
gives a $1$-cocycle of 
$\mathcal O_C(2p)^{\oplus 2}
\simeq N_{S/\mathbb P^4}\big{\vert}_C^{\vee}(2)(2p)
$
with respect to $\mathfrak U_2$.
Now we make a change on open coverings of $C$.
We consider another open affine covering
$$
\mathfrak U_3:=\left\{
D(x_j'), D(\Delta_j) \bigm| j=0,1,2
\right\}
$$
of $C$. Then both $\mathfrak U_1$ and $\mathfrak U_2$ are refinement
of $\mathfrak U_3$. There are isomorphisms between
$$
{\check H}^1(\mathfrak U_k, \mathcal O_C(2p)^{\oplus 2}) \quad (1 \le k \le 3)
$$
induced by natural maps
$$
C^{\bullet}(\mathfrak U_1,\mathcal O_C(2p)^{\oplus 2})
\longleftarrow
C^{\bullet}(\mathfrak U_3,\mathcal O_C(2p)^{\oplus 2})
\longrightarrow
C^{\bullet}(\mathfrak U_2,\mathcal O_C(2p)^{\oplus 2})
$$
of {\v C}ech complexes with respect to $\mathfrak U_k$ ($1\le k\le 3$).
Moreover, by the above computations, 
we see that the 1-cocycle representing $u_i \cup \mathbf e'$
can be taken from the one in  $C^1(\mathfrak U_3,\mathcal O_C(2p)^{\oplus 2})$,
and mapped to $\xi_i(- \hat t_2\big{\vert}_C, \hat t_1\big{\vert}_C)$
in $C^1(\mathfrak U_2,\mathcal O_C(2p)^{\oplus 2})$.
Hence we obtain
$r(u_i \cup \mathbf e',2p)
=\xi_i(-t_2\big{\vert}_C, t_1\big{\vert}_C)
=\xi_i \mathbf t\big{\vert}_C
\in H^1(\mathcal O_C(2p)^{\oplus 2})$ and we have proved (1).
Finally, we prove (2). 
We see by \eqref{eqn:second product} that
$\ob_S(\alpha)=(u_1,u_2)\cup \mathbf e' \cup \alpha'
=(u_1\cup \mathbf e' \cup \alpha', u_2\cup \mathbf e' \cup \alpha')$.
Then since the inclusion $\mathcal O_C^{\oplus 2}(2p)
\hookrightarrow N_{C/\mathbb P^4}^{\vee}(2)(2p)$ 
of sheaves is induced by $\pi_{C/S}^{\vee}$
we have 
$$
r(u_i \cup \mathbf e' \cup \alpha',2p)
=r(u_i \cup \mathbf e',2p)\cup \alpha'
=\pi_{C/S}^{\vee}(\xi_i \mathbf t\big{\vert}_C)\cup \alpha'
=\xi_i \pi_{C/S}^{\vee}(\mathbf t\big{\vert}_C)\cup \alpha'
$$
in $H^1(C,\mathcal O_C(2)(2p))$.
Then there exists a commutative diagram
$$
\begin{CD}
  H^1(C,N_{C/\mathbb P^4}^{\vee}(1)(2p)) @. \times @. H^0(C,N_{C/\mathbb P^4})
  @>{\cup}>> H^1(C,\mathcal O_C(1)(2p)) \\
  @A{\pi_{C/S}^{\vee}}AA @.  @V{\pi_{C/S}}VV \Vert \\
  H^1(C,N_{S/\mathbb P^4}\big{\vert}_C^{\vee}(1)(2p)) @. \times @. 
  H^0(C,N_{S/\mathbb P^4}\big{\vert}_C) @>{\cup}>> H^1(C,\mathcal O_C(1)(2p)).
\end{CD}
$$
Therefore, we obtain 
$$
\pi_{C/S}^{\vee}(\mathbf t\big{\vert}_C) \cup \alpha'
=\mathbf t\big{\vert}_C \cup \pi_{C/S}(\alpha).
$$
Thereby we conclude that
$$
r(\ob_S(\alpha),2p)=\xi \pi_{C/S}^{\vee}(\mathbf t\big{\vert}_C) \cup \alpha'
=\xi \left(\mathbf t\big{\vert}_C \cup \pi_{C/S}(\alpha)\right)
$$
as in the statement.
\qed

\medskip

\paragraph{\bf Proof of Theorem~\ref{thm:obstruction}~(2)}

It suffices for the proof to show 
that the cohomology class $r(\ob_S(\alpha),2p)$
in $H^1(C,N_{S/\mathbb P^4}\big{\vert}_C(2p))$ is nonzero.
We consider the cup product of $r(\ob_S(\alpha),2p)$ with $\mathbf k_C$,
where $\mathbf k_C$ is the extension class of the exact sequence \eqref{ses:CE}.
Then $r(\ob_S(\alpha),2p) \cup \mathbf k_C$
belongs to $H^2(S,N_{S/\mathbb P^4}(2E-C))$
and we have
$$
r(\ob_S(\alpha),2p) \cup \mathbf k_C
=\xi \left(
  \mathbf t\big{\vert}_C \cup \pi_{C/S}(\alpha)
\right)\cup \mathbf k_C 
=\xi \mathbf t \cup 
  \left(
    \pi_{C/S}(\alpha) \cup \mathbf k_C
  \right) 
$$
by \eqref{eqn:exterior component} in Lemma~\ref{lem:key lemma} and 
the associative property on cup products; indeed, there exists 
a commutative diagram
\begin{equation}\label{diag:res-coboundary1}
\begin{array}{ccccc}
  && \pi_{C/S}(\alpha) &&  \\ [-10pt]
  && \rotatebox{-90}{$\in$} &&  \\ [4pt]
  H^1(C,N_{S/\mathbb P^4}\big{\vert}_C^{\vee}(1)(2p)) & \times & H^0(C,N_{S/\mathbb P^4}\big{\vert}_C) 
  & \overset{\cup}{\longrightarrow} & H^1(C,\mathcal O_C(1)(2p))\\
 \mapup{|_C} && \mapdown{\cup \, \mathbf k_C} && 
  \mapdown{\cup \, \mathbf k_C} \\
 H^1(S,N_{S/\mathbb P^4}^{\vee}(1)(2E)) & \times & H^1(S,N_{S/\mathbb P^4}(-C))
  & \overset{\cup}{\longrightarrow} & H^2(S,\mathcal O_S(1)(2E-C)). \\
 \rotatebox{90}{$\in$} && && \\ [-4pt] 
  \mathbf t && && \\ 
\end{array}
\end{equation}
Here we recall again that 
$r(\pi_{C/S}(\alpha),p) \in H^0(C,N_{S/\mathbb P^4}(E)\big{\vert}_C)$
lifts to the global section $\beta$ of $N_{S/\mathbb P^4}(E)$.
Moreover, applying Lemma~\ref{lem:restriction to C and E}~(2)-[iii] to
$\gamma \cup \mathbf k_C$ with $\mathcal F=N_{S/\mathbb P^4}$ and 
$\gamma=\pi_{C/S}(\alpha)$, we obtain
$$
\pi_{C/S}(\alpha) \cup \mathbf k_C = \beta_E' \cup \mathbf k_E,
$$
where $\beta_E'$ is a nonzero global section of the trivial bundle
$N_{S/\mathbb P^4}(E-C)\big{\vert}_E \simeq \mathcal O_E^{\oplus 2}$
on $E$ and such that $r(\beta'_E,C)=\beta\big{\vert}_E$
in $H^0(N_{S/\mathbb P^4}(E)\big{\vert}_E)$.
Then there exists a commutative diagram
\begin{equation}\label{diag:res-coboundary2}
\begin{array}{ccccc}
  && \beta_E' &&  \\ [-10pt]
  && \rotatebox{-90}{$\in$} &&  \\ [4pt]
  H^1(E,N_{S/\mathbb P^4}\big{\vert}_E^{\vee}(1)(2E)) & \times & H^0(E,N_{S/\mathbb P^4}(E-C)\big{\vert}_E) 
  & \overset{\cup}{\longrightarrow} & H^1(E,\mathcal O_E(1)(3E-C))\\
 \mapup{|_E} && \mapdown{\cup \, \mathbf k_E} && 
  \mapdown{\cup \, \mathbf k_E} \\
 H^1(S,N_{S/\mathbb P^4}^{\vee}(1)(2E)) & \times & H^1(S,N_{S/\mathbb P^4}(-C))
  & \overset{\cup}{\longrightarrow} & H^2(S,\mathcal O_S(1)(2E-C)). \\
 \rotatebox{90}{$\in$} && && \\ [-4pt] 
  \mathbf t && && \\ 
\end{array}
\end{equation}
similar to \eqref{diag:res-coboundary1}. Again by commutativity,
we have
\begin{equation}
  \label{eqn:final cup products}
  r(\ob_S(\alpha),2p) \cup \mathbf k_C
  =\xi\mathbf t\cup 
  \left(
    \beta_E' \cup \mathbf k_E
  \right) 
  =\xi\big{\vert}_E\left(
    \mathbf t\big{\vert}_E \cup \beta_E' 
  \right) 
  \cup \mathbf k_E.
\end{equation}
Since $\beta_E'$ is a nonzero global section of a trivial bundle,
Lemma~\ref{lem:linear independence of coh classes}
shows that $\mathbf t\big{\vert}_E \cup \beta_E'\ne 0$.
We recall the restriction $\xi\big{\vert}_E$ of $\xi$ to $E$
has a zero at $p=C\cap E$, which is general in $|\mathcal O_E(1)|$
by Claim~\ref{claim:positivity and vanishing}.
The Serre duality shows that the pairing
$$
H^0(E,\mathcal O_E(1)) \times H^1(E,\mathcal O_E(1)(3E-C))
\overset{\times}{\longrightarrow} H^1(E,K_E)\simeq k
$$
is non-degenerate and hence 
$\mathbf t\big{\vert}_E \cup \beta_E'$ 
multiplied by $\xi\big{\vert}_E$
is nonzero. Therefore, it suffices to show 
that the coboundary map
$$
H^1(E,N_{S/\mathbb P^4}\big{\vert}_E(3E-C))
\overset{\cup \mathbf k_E}{\longrightarrow}
H^2(S,N_{S/\mathbb P^4}(2E-C))
$$
is injective. Indeed, since $\chi(-2K_S-C)\ge 0$,
it follows from \S\ref{subsec:delpezzo}({\bf P6}) that
$$
H^1(S,N_{S/\mathbb P^4}(3E-C))\simeq H^1(S,-2K_S+3E-C)^{\oplus 2}=0.
$$
Then by \eqref{eqn:final cup products}, 
we have proved that $\ob_S(\alpha)\ne 0$
and hence $\ob(\alpha)\ne 0$.
\qed

\subsection{Non-reduced components}
\label{subsec:non-reduced}

In this section, we prove Theorem~\ref{thm:main2}.
Let $C$ be a smooth connected curve of degree $d$ and genus $g$
contained in a smooth complete intersection $S=S_{2,2}$ in $\mathbb P^4$.

\begin{lem}
  \label{lem:g>=2d-12}
  Let $L$ denote the class of $\mathcal O_S(C+2K_S)$ in $\Pic S$ 
  as in \S~\ref{subsec:obstructions}. Then
  \begin{enumerate}
    \item $\chi(S,-L)=g-2d+12$.
    \item Suppose that $d>8$ and $g \ge 2d-12$.
    \begin{enumerate}
      \item If $L$ is nef,
      then $C$ is $2$-normal.
      \item If $C$ is not $2$-normal,
      then $L+K_S$ is effective and 
      $$k:=h^1(\mathcal I_C(2))=h^0(F,\mathcal O_F),$$
      where $F$ is the fixed part of $|L|$.
      If moreover $C$ is linearly normal then 
      there exist disjoint $k$ lines $E_i$ ($1 \le i \le k$) on $S$
      such that $F=E_1+\dots+E_k$.
    \end{enumerate}
  \end{enumerate}
\end{lem}
\Proof
The proof is similar to that of \cite[Lemma~3.1]{Nasu9}.
By Riemann-Roch theorem on $S$, 
$$
\chi(S,-L)=(C+2K_S)(C+3K_S)/2+1=g-2d+12
$$
and we obtain (1). Suppose that $d>8$ and $g \ge 2d-12$.
Then $H^0(S,-L)=0$ by Lemma~\ref{lem:uniqueness of pencil} and $\chi(S,-L)\ge 0$ by (1).
We recall that 
$H^1(\mathcal I_C(2))\simeq H^1(S,-L)$ by \eqref{isom:normality}.
Thereby [i] follows from \S\ref{subsec:delpezzo}({\bf P4}). We prove [ii].
Suppose that $C$ is not $2$-normal, i.e.,~$H^1(S,-L)\ne 0$.
Then since $\chi(S,-L)\ge 0$, we obtain $H^2(S,-L)\ne 0$.
Therefore $L+K_S$ is effective by Serre duality and so is $L$.
Since $L$ is not nef by [i], the complete linear system
$|L|$ has the nonzero fixed part
$F=-\sum_{\ell} (L.\ell) \ell$, which is a sum of lines $\ell$
on $S$ such that $L.\ell<0$.
Then $L$ is big by \S\ref{subsec:delpezzo}({\bf P6}) and hence
$k=h^1(S,-L)=h^0(F,\mathcal O_F)$ by \S\ref{subsec:delpezzo}({\bf P4}).
Finally suppose that $C$ is linearly normal.
Then $H^1(S,-C-K_S) \simeq H^1(\mathcal I_C(1))=0$ by \eqref{isom:normality}, 
while $C+K_S$ is effective and big by $C+K_S=L-K_S$.
Then this implies that $C+K_S$ is nef, because otherwise
$H^1(S,-C-K_S) \simeq H^0(F',\mathcal O_{F'})\ne 0$ with $F'=\Bs|C+K_S|$
by \S\ref{subsec:delpezzo}({\bf P4}). Since $L=C+2K_S$ is not nef, 
$F$ is a sum of lines $\ell$ on $S$ such that $L.\ell=-1$
and $k$ ($=h^0(F,\mathcal O_F)$) coincides with the number of lines in $F$.
\qed

\medskip

\paragraph{\bf Proof of Theorem~\ref{thm:main2}}
Let $C$ and $S$ be as in the theorem. 
Then by Lemma~\ref{lem:g>=2d-12}, we have
$L+K_S \ge 0$ and 
$h^1(\mathcal I_C(2))=h^0(F,\mathcal O_F)$, where $F$
is the fixed part of $|L|$. 
Since $h^1(\mathcal I_C(2))=1$,
there exists a line $E$ on $S$ such that $F=E$.
Here we note that $F$ does not contain any double line $2E$,
because otherwise $h^0(F,\mathcal O_F) \ge h^0(2E,\mathcal O_{2E})=3$.
We now apply Theorem~\ref{thm:obstruction} to $C$ and $S$.
The first statement of the theorem shows that
the tangent map $p_1$ of the first projection 
$pr_1: \HF^{sc} \mathbb P^4 
\rightarrow \Hilb^{sc} \mathbb P^4$ at $(C,S)$ 
is not surjective and its cokernel is of dimension $2$.
Let $C'$ be a general member of $W$.
Then $C'$ is general also
in the class $[C']$ of $C'$ in $\Pic S$ by dimensions.
Then it follows from the second statement of Theorem~\ref{thm:obstruction} that
if $\alpha \in H^0(N_{C'/\mathbb P^4})$ is not 
contained in the image of $p_1$,
then its primary obstruction $\ob(\alpha)$ is nonzero in 
$H^1(N_{C'/\mathbb P^4})$.
This implies that 
$\Hilb^{sc} \mathbb P^4$ is singular at $[C']$ by infinitesimal lifting property.
It follows from Lemma~\ref{lem:maximality} that
the $S$-maximal family $W_{C,S}$ ($=W_{C',S'}\supset W$)
is an irreducible component of
$(\Hilb^{sc} \mathbb P^4)_{\red}$.
Since $W$ is dense in $W_{C,S}$ by Remark~\ref{rmk:density and finiteness},
$\Hilb^{sc} \mathbb P^4$ is generically singular 
and hence generically non-reduced along $W_{C,S}=\overline W$. \qed

\medskip

We give a sufficient condition for a $S_{2,2}$-maximal family
to be an irreducible component of $(\Hilb^{sc} \mathbb P^4)_{\red}$ 
and discuss the generic smoothness of $\Hilb^{sc} \mathbb P^4$ along 
the component, in terms of the corresponding 
$6$-tuples $(a;b_1,\dots,b_5)$ of integers
(cf.~Lemma~\ref{lem:maximal families and 6-tuples}).
We refer to \cite{Kleppe87,Ellia87,Nasu1,Kleppe-Ottem15,Kleppe17,Nasu9} 
for corresponding results on 
$S_3$-maximal families in $\Hilb^{sc} \mathbb P^3$.

\begin{thm}
  \label{thm:conclusion}
  Let $W:=W(a;b_1,\dots,b_5) \subset \Hilb_{d,g}^{sc} \mathbb P^4$ 
  be a $S_{2,2}$-maximal family of smooth connected curves 
  of degree $d$ and genus $g$ in $\mathbb P^4$ whose general member $C$
  is contained in $S_{2,2} \subset \mathbb P^4$ and
  has the standard coordinate $(a;b_1,\dots,b_5)$ 
  in $\Pic S_{2,2}\simeq \mathbb Z^6$
  (see Definition~\ref{dfn:standard coordinates}).
  Suppose that $d>10$ and $g \ge 2d-12$. Then
  \begin{enumerate}
    \item If $b_5\ge 2$, then 
    $W$ is an irreducible component of $\Hilb_{d,g}^{sc} \mathbb P^4$
    and $\Hilb_{d,g}^{sc} \mathbb P^4$ is generically smooth along $W$.
    \item If $b_5=1$ and $b_4\ge 2$, 
    then $W$ is an irreducible component of $(\Hilb_{d,g}^{sc} \mathbb P^4)_{\red}$
    and $\Hilb_{d,g}^{sc} \mathbb P^4$ is generically non-reduced along $W$.
    \item If $b_5=0$, then $W$
    is not an irreducible component $(\Hilb_{d,g}^{sc} \mathbb P^4)_{\red}$,~i.e.,~
    there exists an irreducible component of
    $(\Hilb_{d,g}^{sc} \mathbb P^4)_{\red}$ containing $W$ strictly.
  \end{enumerate}
\end{thm}
\Proof
Put $S=S_{2,2}$ and $L:=[\mathcal O_S(C+2K_S)]$ in $\Pic S$.
If $b_5 \ge 2$ then $L$ is nef, and hence $C$ is $2$-normal by Lemma~\ref{lem:g>=2d-12}.
Thus (1) follows from Theorem~\ref{thm:main1}.
If $b_5=1$ and $b_4\ge 2$, 
then the fixed part of $|L|$ is equal to a line $E$ ($=\mathbf e_5$) on $S$.
We see that $L^2=2g-3d+14\ge 2(2d-12)-3d+14=d-10>0$.
Then $h^1(\mathcal I_C(2))=h^1(S,-L)=h^0(E,\mathcal O_E)=1$.
Thus (2) follows from Theorem~\ref{thm:main2}.
Finally, if $b_5=0$, then $C$ is not linearly normal by $C+K_S \ge 0$ and $(C+K_S)^2>0$,
both of which easily follow from Riemann-Roch theorem on $S$.
Then by Lemma~\ref{lem:not 1-normal}, for every general member $C$ of $W$,
there exists a global deformation of $C'$ of $C$ in $\mathbb P^4$ such that
$C'$ is not contained in any smooth complete intersection $S_{2,2}$ in $\mathbb P^4$.
Thus (3) follows.
\qed

\medskip

The next example shows that there exist
infinitely many generically non-reduced components
of $\Hilb^{sc}\mathbb P^4$, whose general member is contained in 
a smooth complete intersection $S_{2,2} \subset \mathbb P^4$.

\begin{ex}
  \label{ex:non-reduced comps.inf.many}
  Let $\lambda \ge 0$ be any integer and let
  $$
  W_\lambda:=W(\lambda+9;\lambda+2,2,2,2,1) 
  \subset \Hilb_{d,g}^{sc}\mathbb P^4.
  $$
  Then $d=2\lambda+18$ and $g=6\lambda+24$,
  and thereby $\dim W_\lambda=d+g+25=8\lambda+67$. 
  Since $b_5=1$, $b_4 \ge 2$
  and $g-2d+12=2\lambda \ge 0$, we see by Theorem~\ref{thm:conclusion}
  that $W_\lambda$ is an irreducible component of
  $(\Hilb^{sc}\mathbb P^4)_{\red}$ and
  $\Hilb^{sc}\mathbb P^4$ is generically non-reduced along $W_\lambda$.
\end{ex}

Table~\ref{table:list of maximal families} contains any solutions of
\eqref{standard-prescribed-nef-big}
when $(d,g)$ belongs to the region (II) of Figure~\ref{fig:d-g-region}
and $14 \le d \le 18$.
The corresponding $S$-maximal family $W(a;b_1,\dots,b_5)$
becomes either a generically smooth component,
a generically non-reduced component,
or a proper closed subset of an irreducible component 
of the Hilbert scheme $\Hilb_{d,g}^{sc} \mathbb P^4$ 
by Theorem~\ref{thm:conclusion} as in the table.

\begin{table}
  \caption{$S_{2,2}$-maximal families in $\Hilb^{sc}_{d,g} \mathbb P^4$}
  \begin{center}
    \begin{tabular}{|c|c|l|}
      \hline
      $(d,g)$ & $(a;b_1,b_2,b_3,b_4,b_5)$ & $W(a;b_1,b_2,b_3,b_4,b_5)$\\
      \hline
      $(14,16)$ & $(8;2,2,2,2,2)$ & gen.~smooth component \\
      $(14,16)$ & $(9;4,3,2,2,2)$ & gen.~smooth component \\
      $(14,16)$ & $(9;3,3,3,3,1)$ & gen.~non-reduced component \\
      \hline
      $(15,18)$ & $(9;4,2,2,2,2)$ & gen.~smooth component \\
      $(15,18)$ & $(9;3,3,3,2,1)$ & gen.~non-reduced component \\
      \hline
      $(16,20)$ & $(9;3,3,2,2,1)$ & gen.~non-reduced component \\
      $(16,20)$ & $(10;5,3,2,2,2)$ & gen.~smooth component \\
      \hline
      $(16,21)$ & $(9;3,2,2,2,2)$ & gen.~smooth component \\
      $(16,21)$ & $(10;4,4,2,2,2)$ & gen.~smooth component \\
      $(16,21)$ & $(10;4,3,3,3,1)$ & gen.~non-reduced component \\
      \hline
      $(17,22)$ & $(9;3,2,2,2,1)$ & gen.~non-reduced component \\
      $(17,22)$ & $(10;5,2,2,2,2)$ & gen.~smooth component \\
      $(17,22)$ & $(10;4,4,2,2,1)$ & gen.~non-reduced component \\
      \hline
      $(17,23)$ & $(9;2,2,2,2,2)$ & gen.~smooth component \\
      $(17,23)$ & $(10;4,3,3,2,1)$ & gen.~non-reduced component \\
      \hline
      $(17,24)$ & $(10;4,3,2,2,2)$ & gen.~smooth component \\
      $(17,24)$ & $(10;3,3,3,3,1)$ & gen.~non-reduced component \\
      \hline
      $(18,24)$ & $(9;2,2,2,2,1)$ & gen.~non-reduced component \\
      $(18,24)$ & $(10;4,3,3,1,1)$ & unknown \\
      $(18,24)$ & $(10;3,3,3,3,0)$ & non-component \\
      $(18,24)$ & $(11;6,3,2,2,2)$ & gen.~smooth component \\
      \hline
      $(18,25)$ & $(10;4,3,2,2,1)$ & gen.~non-reduced component \\
      \hline
      $(18,26)$ & $(10;4,2,2,2,2)$ & gen.~smooth component \\
      $(18,26)$ & $(10;3,3,3,2,1)$ & gen.~non-reduced component \\
      $(18,26)$ & $(11;5,4,2,2,2)$ & gen.~smooth component \\
      $(18,26)$ & $(11;5,3,3,3,1)$ & gen.~non-reduced component \\
      \hline
      $(18,27)$ & $(10;3,3,2,2,2)$ & gen.~smooth component \\
      $(18,27)$ & $(11;5,3,3,2,2)$ & gen.~smooth component \\
      $(18,27)$ & $(11;4,4,3,3,1)$ & gen.~non-reduced component \\
      \hline
    \end{tabular}
    \label{table:list of maximal families}
  \end{center}
\end{table}

\bibliography{mybib}
\bibliographystyle{abbrv}

\end{document}